\documentclass[12pt]{article}
\usepackage{latexsym,amssymb,amsmath}
\usepackage{amsthm, amstext}
\usepackage{array, amsfonts, mathrsfs}

\usepackage[dvips]{graphicx,color}

\newtheorem{Theorem}{Theorem}
\newtheorem{Lemma}{Lemma}

\newtheorem{Proposition}{Proposition}
\def\blue{\color{blue}}

\numberwithin{equation}{section}

\date{}
\begin{document}

\author{M.I.Belishev\thanks {Saint-Petersburg Department of
                 the Steklov Mathematical Institute, belishev@pdmi.ras.ru;
                 Saint-Petersburg State University, 7/9 Universitetskaya
                 nab., St. Petersburg, 199034 Russia, m.belishev@spbu.ru.
                 Supported by grants RFBR 14-01-00535А and SPbGU
                 6.38.670.2013.}\, and A.F.Vakulenko\thanks{Saint-Petersburg Department of
                 the Steklov Mathematical Institute, vak@pdmi.ras.ru}}

\title{On Characterization of Inverse Data in the Boundary Control Method}
\maketitle

\begin{abstract}
We deal with a dynamical system
 \begin{align*}
& u_{tt}-\Delta u+qu=0 && {\rm in}\,\,\,\Omega \times (0,T)\\
& u\big|_{t=0}=u_t\big|_{t=0}=0 && {\rm in}\,\,\,\overline
\Omega\\
& \partial_\nu u = f && {\rm in}\,\,\,\partial\Omega \times
[0,T]\,,
 \end{align*}
where $\Omega \subset {\mathbb R}^n$ is a bounded domain, $q \in
L_\infty(\Omega)$ a real-valued function, $\nu$ the outward normal
to $\partial \Omega$, $u=u^f(x,t)$ a solution. The input/output
correspondence is realized by a response operator $R^T: f \mapsto
u^f\big|_{\partial\Omega \times [0,T]}$ and its relevant extension
by hyperbolicity $R^{2T}$. Ope\-rator $R^{2T}$ is determined by
$q\big|_{\Omega^T}$, where $\Omega^T:=\{x \in \Omega\,|\,\,{\rm
dist\,}(x,\partial \Omega)<T\}$. The inverse problem is: Given
$R^{2T}$ to recover $q$ in $\Omega^T$. We solve this problem by the
boundary control method and describe the {\it ne\-ces\-sary and
sufficient} conditions on $R^{2T}$, which provide its solvability.
\end{abstract}

\section{Introduction}\label{sec Introduction}
\subsubsection*{Motivation}
The problem, which the paper is devoted to, was solved about 20
years ago by the BC-method, which is an approach to inverse
problems (IPs) based on their relations to control and system
theory \cite{BKach94, BIP97, BIP07}. However, in IP-community,
there are a few versions of what 'to solve an inverse problem'
means. The versions may be ordered by levels as follows:

$1.$\,\,\,to establish injectivity of the correspondence
`parameters under reconstruction $\to$ inverse data', what allows
one to claim that the data {\it determine} the parameters

$2.$\,\,\,to elaborate an efficient (preferably, realizable
numerically) {\it procedure}, which determines the parameters from
the data \footnote{surely, we mean the mathematically rigorous
approaches)}

$3.$\,\,\,to provide a {\it data characterization}, i.e., describe
the necessary and sufficient conditions on the data, which ensure
solvability of the given inverse problem.

\noindent Typically, $\{i+1\}$-th level is stronger and richer in
content than $i$-th one. Respectively, to reach the next level
(especially, in multidimensional IPs) is more difficult. The
BC-method firmly keeps level 2 (see \cite{BIP07, BIKS}). In the
mean time, it provides a data characterization in important
one-dimensional problems: see \cite{BMikh, BPest}.

Regarding level 3 in {\it multidimensional} IPs, there is
substantial gap between the frequency-domain and time-domain
problems. In the first ones, the results on the data
characterization are much more promoted and successful (see
\cite{Fadd, HN, Newt, Nov} and other). In time-domain problems,
such results also do exist (see, e.g., \cite{Romanov}) but are not
so deep and systematic. Our paper is an attempt to reduce the
above-mentioned gap by the use of the BC-method.

\subsubsection*{Contents and results}
$\bullet$\,\,\,We develop a general approach proposed in
\cite{DSBC} and apply it to a concrete time-domain inverse problem
for the wave equation with a potential. The approach elaborates
the well-known and deep relations between inverse problems and
triangular factorization of operators in the Hilbert space
\cite{Fadd, BIP97, DSBC, BPush}.
\smallskip

\noindent$\bullet$\,\,\,In sections \ref{sec Geometry} and
\ref{sec Dynamics}, a forward problem is considered. With the
problem one associates a relevant dynamical system. The system is
endowed with standard control theory attributes: spaces and
operators. In particular, a so-called {\it extended response
operator} $R^{2T}$ is introduced. It realizes the input/state
correspondence and later on plays a role of the data in the
inverse problem. The key property of the system is a {\it local
boundary controllability}, which is relayed upon the fundamental
Holmgren-John-Tataru uniqueness theorem \cite{Tat}. It plays a
crucial role in all versions of the BC-method.

Geometrical Optics (GO) describes propagation of wave field jumps
in the system. A noticeable fact is that the GO-formulas are well
interpreted in operator theory terms: they provide existence of a
{\it diagonal} of the control operator and time derivative
composition.
\smallskip

\noindent$\bullet$\,\,\,In section \ref{sec Visualization},  we
present a BC-procedure, which recovers the potential from the
given $R^{2T}$. Then we prove Theorem \ref{T1}, which is the main
result. It provides a list of necessary and sufficient conditions
on an operator ${\cal R}^{2T}$ to be an extended response
operator.

The necessity is simple: the proof just summarizes the properties
of $R^{2T}$ stated in the forward problem. The sufficiency is
richer in content. The proof is constructive: we start with an
operator ${\cal R}^{2T}$ obeying all the conditions, and construct
a system with the response operator $R^{2T}={\cal R}^{2T}$. In
construction we follow the BC-procedure, which solves the IP.

In conclusion (section \ref{Comments}), a self-critical discussion
of the obtained results is provided.

\subsubsection*{Acknowledgements}
We wold like to thank R.G.Novikov and V.G.Romanov for kind
consultations.

\section{Geometry}\label{sec Geometry}
All the functions, function classes and spaces are {\it real}.

\subsubsection*{Domain and subdomains} Let ${\Omega}\subset {\mathbb
R}^n$ be a bounded domain with the boundary $\Gamma\in C^\infty$.
By ${\rm d}(a,b)$ we denote an {\it intrinsic distance} in
$\Omega$, which is defined via the length of smooth curves lying
in $\overline \Omega$ and connecting $a$ with $b$.

For a subset $A \subset \overline{\Omega}$, we denote its metric
neighborhoods by
$$
{\Omega}^r_A := \{x \in {\Omega}\,|\,\,{\rm d}(x,A) < r\}, \qquad
r>0.
$$
For $A=\Gamma$, we set ${\Omega}^r:={\Omega}^r_\Gamma$. Later on,
in dynamics, the value
$$
T_*:=\max_{\Omega} \tau(\cdot)=\inf
\{r>0\,|\,\,{\Omega}^r={\Omega} \}
$$
is interpreted as a time needed for the waves moving from $\Gamma$
with the unit speed to fill ${\Omega}$.

A function $\tau(\cdot):={\rm d}(\cdot,\Gamma)$ on $\overline
{\Omega}$ is called an {\it eikonal}. By the definitions, we have
${\Omega}^r=\{x \in {\Omega}\,|~\tau(x)<r\}$. In dynamics, the
eikonal level sets
$$
\Gamma^s:=\{x\in {\Omega}\,|\,\,\tau(x)=s\}, \qquad s \geqslant 0
$$
play the role of the forward fronts of waves moving from $\Gamma$.

\subsubsection*{Semi-geodesic coordinates}
$\bullet$\,\,\,Here we introduce a separation set (cut locus) of
${\Omega}$ with respect to $\Gamma$ (see, e.g, \cite{GKM}) and use
one of its equivalent definitions \cite{Kling}.

A point in $\Omega$ is said to be multiple if it is connected with
$\Gamma$ through more than one shortest geodesics (straight lines
in ${\mathbb R}^n$). Denote by $c_0$ the set of multiple points
and define
$$c\,:=\,\overline c_0.$$
The set $c$ is called a {\it cut locus}. It is 'small':
 \begin{equation}\label{vol c=0}
{\rm vol\,}c\,=\,0\,,
 \end{equation}
and separated from the boundary:
 $$
0\,<\,T_c\,:={\rm d}(c, \Gamma)\,\leqslant T_*\,.
 $$

In addition, note that $\Gamma^s \backslash c$ is a smooth (may
be, disconnected) hyper-surface in ${\Omega}$. If $s<T_c$ then
$\Gamma^s$ is smooth and diffeomorphic to $\Gamma$.
\smallskip

\noindent$\bullet$\,\,\, For any $x \in \overline{\Omega}\,
\backslash c$, there is a unique point $\gamma(x)\in\Gamma$
nearest to $x$. For such an $x$, a pair $(\gamma(x),\tau(x))$
determines its position in $\Omega$ and is said to be the {\it
semi-geodesic coordinates} (sgc). By $x(\gamma, \tau)$ we denote a
point in $\overline{\Omega}\, \backslash c$ with the given sgc
$(\gamma, \tau)$.

In sgc, ${\mathbb R}^n$-volume element in $\Omega$ takes the
well-known form
 \begin{equation}\label{dx=beta dGamma dtau}
dx\,=\,\beta(\gamma, \tau)\,d\Gamma d\tau\,,
 \end{equation}
where $d\Gamma$ is Euclidean surface element on the boundary.
Factor $\beta$ is a Jacobian of the passage from Cartesian
coordinates to sgc.
\smallskip

\noindent$\bullet$\,\,\,Denote $\Sigma^T:=\Gamma \times [0,T)$. A
set
 $$
\Theta\,:=\,{\{(\gamma(x),\tau(x))\,|\,\,x \in \left[\Omega \cup
\Gamma\right]\backslash c\}}\, \subset\, \Sigma^{\,T_*}
 $$
is called a {\it pattern} of ${\Omega}$. Also, we use its parts
 $$
\Theta^T:=\left\{\left(\gamma(x),\tau(x)\right)\,|\,\,x \in
\left[{\Omega^T}\cup\Gamma\right] \backslash c\right\} = \Theta
\cap \Sigma^T, \quad T>0\,.
 $$
For $T<T_c$, one has $\Theta^T = \Sigma^T$.

\subsubsection*{Images}
Fix a positive $T \leqslant T_*$; let $y$ be a function on $\Omega^T
\cup \Gamma$. A function on $\Sigma^T$ of the form
 $$
\tilde y^T(\gamma, \tau)\,:=\,
  \begin{cases}
\beta^{\frac{1}{2}}(\gamma,\tau)y\left(x(\gamma,\tau)\right), &
(\gamma,\tau) \in \Theta^T\\
0, & (\gamma,\tau) \in \Sigma^T \backslash \Theta^T
  \end{cases}
 $$
is said to be an {\it image} of $y$. So, up to the factor
$\beta^{\frac{1}{2}}$, image is just a function written in sgc.

An {\it image operator} $ I^T:\,L_2(\Omega^T) \to
L_2(\Sigma^T),\,\,\,I^Ty\,:=\,\tilde y^T$ is isometric. Indeed,
for $y,v \in L_2(\Omega^T)$ one has
 \begin{align*}
&
\left(y,v\right)_{L_2(\Omega^T)}=\int_{\Omega^T}y(x)\,v(x)\,dx\overset{(\ref{vol
c=0}), (\ref{dx=beta dGamma
dtau})}=\int_{\Theta^T}y(x(\gamma,\tau))\,v(x(\gamma,\tau))\,\beta(\gamma,\tau)d\Gamma\,d\tau=\\
& =\,\left(\tilde y^T, \tilde
v^T\right)_{L_2(\Sigma^T)}\,=\,\left(I^T y, I^T
v\right)_{L_2(\Sigma^T)}\,.
 \end{align*}

As an isometry, $I^T$ obeys ${\rm Ran\,} I^T=\{g \in
L_2(\Sigma^T)\,|\,\,{\rm supp\,} g \subset \overline{\Theta^T}\}$ and
 \begin{equation}\label{I^T properties}
\left(I^T\right)^* I^T\,=\,{\mathbb I}_{L_2(\Omega^T)}\,, \quad
I^T \left(I^T\right)^*\,=\,G_{\Theta^T}\,,
 \end{equation}
where $G_{\Theta^T}$ cuts off functions in $\Sigma^T$ onto
$\Theta^T$.

\section{Dynamics}\label{sec Dynamics}
\subsection{IBV-problem}\label{subsec IBV-problem}
By $\partial_\nu$ we denote a derivative with respect to outward
normal at the boundary $\Gamma$. $H^s(\dots)$ are the standard
Sobolev spaces.
\smallskip

Consider an initial boundary-value problem
 \begin{align}
\label{A1}
 & u_{tt}-\Delta u+qu=0 &&{\rm in}\,\,\,{\Omega} \times
(0,T)\\
\label{A2}
 &  u\big|_{t=0}=u_t\big|_{t=0}=0 \quad &&{\rm in}\,\,\,\overline \Omega\\
\label{A3} &  \partial_\nu u=f \quad &&{\rm
on}\,\,\,\overline{\Sigma^T}\,,
\end{align}
where $q \in L_\infty(\Omega)$ is a function ({\it potential}),
$f$ is a Neumann {\it boundary control}, $u=u^f(x,t)$ is a
solution ({\it wave}). It is a well-posed problem; its solution
possesses the following properties.
\smallskip

\noindent$\bullet$\,\,{\it Regularity.}\,\,\,The map $f \mapsto
u^f$ is continuous from $L_2(\Sigma^T)$ to $C([0,T];
H^{\frac{3}{5}-\varepsilon}(\Omega))$, whereas $f \mapsto
u^f\big|_{\Sigma^T}$ acts continuously from $L_2(\Sigma^T)$ to
$H^{\frac{1}{5}-2\varepsilon}(\Sigma^T)$ \, ($\forall
\varepsilon>0$). Introduce a `smooth' class of controls
 $$
{\cal M}^T\,:=\,\left\{f \in H^2(\Sigma^T)\,|\,\,\rm supp\, f \subset
\Gamma \times
(0,T]\right\}
 $$
and note that each $f \in {\cal M}^T$ vanishes near $t=0$. For $f
\in {\cal M}^T$ one has $u^f \in H^2(\Omega \times [0,T])$. These
facts are taken from \cite{LTr-2} (Theorem A).
\smallskip

\noindent$\bullet$\,\,{\it Locality.}\,\,\,For the hyperbolic
equation (\ref{A1}), the finiteness of the domain of influence
principle holds and implies the following.

Let $\sigma \subset \Gamma$ be an open set. Take a control acting
from $\sigma$, i.e., provided ${\rm supp\,} f \subset \overline \sigma
\times [0,T]$. Then the relation
 \begin{equation}\label{supp u-1}
{\rm supp\,} u^f(\cdot, t)\,\subset\,\overline {\Omega^t_\sigma},
\qquad t  \geqslant 0
 \end{equation}
holds and shows that the waves propagate with the unit speed and
fill the proper metric neighborhood of $\sigma$ in $\Omega$.

By the latter, solution $u^f$ depends on the potential locally
that enables one to restate the problem (\ref{A1})--(\ref{A3}) as
follows:
\begin{align}
\label{B1}
 & u_{tt}-\Delta u+qu=0 &&{\rm in}\,\,\,{\Omega^T} \times
(0,T)\\
\label{B2}
 & u\big|_{t<\tau(x)}\,=\,0 &&{\rm in}\,\,\,\overline{\Omega^T} \times
[0,T]\\
\label{B3} &  \partial_\nu u=f \quad &&{\rm
on}\,\,\,\overline{\Sigma^T}\,.
\end{align}
Such a form emphasizes that $u^f$ is determined by behavior of
potential $q$ in $\Omega^T$ only (does not depend on
$q\big|_{\Omega \backslash \Omega^T}$) that enables one to analyze
wave propagation without leaving $\Omega^T$.
\smallskip

\noindent$\bullet$\,\,{\it Steady-state property.}\,\,\,Introduce
a {\it delay operator} ${\cal T}^T_{T-\xi}$ acting on controls by
the rule
 $$
\left({\cal T}^T_{T-\xi}f\right)(\cdot, t)\,:=\,
  \begin{cases}
0\,, &
0 \leqslant t <T-\xi\\
f(\cdot, t-(T-\xi))\,, & T-\xi \leqslant t \leqslant T
  \end{cases}\,\quad 0 \leqslant t \leqslant T\,.
 $$
Since the operator $-\Delta+q$, which governs the evolution of
waves, does not depend on time, one has
 \begin{align}
\notag & u^{{\cal T}^T_{T-\xi}f}(\cdot, T)=u^f(\cdot, \xi)\,,
\quad 0 \leqslant \xi \leqslant T\,; \\
\label{Delay rel}
 & u^{f_t}\,=\,u^f_t\,,\,\,\,
u^{f_{tt}}\,=\,u^f_{tt}\,\overset{(\ref{A1})}=\left(\Delta-q\right)u^f
\qquad {\rm for\,\,} f \in {\cal M}^T\,,
 \end{align}
where the first relation implies the others.

\subsection{System $\alpha^T$}\label{subsec System alpha^T}
Here we consider problem (\ref{B1})--(\ref{B3}) as a dynamical
system, name it by $\alpha^T$, and endow with standard attributes
of control and system theory: spaces and operators.

\subsubsection*{Spaces and subspaces}
A space of controls ${\cal F}^T:=L_2(\Sigma^T)$ is called an {\it
outer space} of the system. It contains an increasing family of
subspaces, which consist of the delayed controls:
 $$
{\cal F}^{T,\,\xi}\,:=\,\left\{f \in {\cal F}^T\,|\,\,{\rm supp\,} f\subset \Gamma \times
[T-\xi,T]\right\}\,=\,{\cal T}^T_{T-\xi}{\cal F}^T, \qquad 0 \leqslant \xi \leqslant
T\,.
 $$
With an open $\sigma \subset \Gamma$ one associates the subspaces
of controls
 $$
{\cal F}^{T,\,\xi}_\sigma\,:=\,\left\{f \in {\cal F}^T\,|\,\,{\rm supp\,} f\subset
\overline\sigma \times [T-\xi,T]\right\}\,, \qquad 0 \leqslant \xi \leqslant
T\,,
 $$
which act from $\sigma$.

A space ${\cal H}^T=L_2(\Omega^T)$ is said to be {\it inner};
waves $u^f(\cdot,t)$ are regarded as its elements ({\it states})
depending on time. It contains an increasing family of subspaces
 $$
{\cal H}^\xi\,:=\,\{y \in {\cal H}^T\,|\,\, {\rm supp\,} y \subset
\overline{\Omega^T}\,\}\,, \qquad 0 \leqslant \xi \leqslant T\,.
 $$
Also, with $\sigma \subset \Gamma$ we associate the subspaces
 $$
{\cal H}^\xi_\sigma\,:=\,\{y \in {\cal H}^T\,|\,\, {\rm supp\,} y \subset
\overline{\Omega^T_\sigma}\,\}\,, \qquad 0 \leqslant \xi \leqslant T\,.
 $$

By locality property (\ref{supp u-1}) and the first relation in
(\ref{Delay rel}), if $f \in {\cal F}^{T,\,\xi}_\sigma $ then $u^f(\cdot,T)\in {\cal H}^\xi_\sigma$.

\subsubsection*{Control operator}
$\bullet$\,\,\,In system $\alpha^T$, an input/state correspondence
is realized by a {\it control operator} $W^T: {\cal F}^T \to {\cal
H}^T$
 $$
W^T f\,:=\,u^f(\cdot,T)\,.
 $$
By the above mentioned regularity properties of solutions to
(\ref{A1})--(\ref{A3}), it acts continuously from ${\cal F}^T$ to
$H^{\frac{3}{5}-\varepsilon}(\Omega)$. Hence, for any $T>0$, $W^T$
is a compact operator.
 \begin{Lemma}\label{L1}
For  $T<T_*$, the control operator  is injective: ${\rm Ker\,}
W^T=\{0\}$.
 \end{Lemma}
\noindent$\square$\,\,\,\,Let $T<T_*$, so that $\Omega \setminus
\overline{\Omega^T}$ is an open set. Let $f \in {\rm Ker\,} W^T=\{0\}$,
so that $u^f(\cdot,T)=0$. Define a function $U$ in $\Omega \times
{\mathbb R}$ by
 $$
U(\cdot,t)\,:=\,
       \begin{cases}
 0\,, & -\infty<t<0\\
 u^f(\cdot,t)\,, & 0 \leqslant t \leqslant T\\
 - u^f(\cdot,2T-t)\,, & T \leqslant t \leqslant 2T\\
 0\,, & -\infty<t<0\,.
       \end{cases}
 $$
Owing to $u^f(\cdot,T)=0$, such an extension of $u^f$ does not
violate its regularity. As a consequence, the extension satisfies
 \begin{equation*}
 U_{tt}-\Delta U + q U=0 \quad {\rm in}\,\,\,\Omega \times
 \mathbb R\,, \qquad U(\cdot, t)\big|_{\Omega \setminus
 \Omega^T}=0\,.
 \end{equation*}
 Applying the Fourier transform $U(\cdot,t)\mapsto \check
 U(\cdot,\omega)$, we get
\begin{equation*}
 -\omega^2 \check U-\Delta \check U + q \check U=0 \quad {\rm in}\,\,\,\Omega\,,
 \qquad \check U(\cdot, \omega)\big|_{\Omega \setminus
 \Omega^T}=0\,.
 \end{equation*}
Thus, for any $\omega \in \mathbb R$, $\check U(\cdot,\omega)$
satisfies an {\it elliptic} equation and vanishes on an {\it open}
set. By the well-known uniqueness theorem, the latter implies
$\check U(\cdot,\omega)=0$ everywhere in $\Omega$. Returning to
the Fourier original, we get $U(\cdot,t)=0$ for all $t$ and arrive
at $f=\partial_\nu u^f\big|_{\Sigma^T}=\partial_\nu
U\big|_{\Sigma^T}=0$. Thus, $f \in {\rm Ker\,} W^T$ implies $f=0$. \,\,
$\blacksquare$
\smallskip

\noindent$\bullet$\,\,\, The locality property (\ref{supp u-1})
and delay relation (\ref{Delay rel}) lead to the embedding
 \begin{equation}\label{W^T FTxis subset Hxis}
 W^T {\cal F}^{T,\,\xi}_\sigma\,\subset \, {\cal H}^\xi_\sigma\,, \qquad 0 \leqslant \xi \leqslant T\,,
 \end{equation}
which is just a consequence of the finiteness of the wave
propagation speed. The fact, which plays a crucial role in the
BC-method, is that this embedding is {\it dense}: the relation
 \begin{equation}\label{W^T FTxis = Hxis}
 \overline{W^T {\cal F}^{T,\,\xi}_\sigma}\,= \,{\cal H}^\xi_\sigma\,, \qquad 0 \leqslant \xi \leqslant T
 \end{equation}
is valid for any $T>0$ and open $\sigma \subseteq \Gamma$. In
control theory this fact is referred to as a {\it local
approximate boundary controllability} of system $\alpha^T$; it is
derived from the fundamental Holmgren-John-Tataru uniqueness
theorem \cite{BIP97, Tat}.
\smallskip

\noindent$\bullet$\,\,\,The following fact will be required in the
data characterization. A multiplication of functions by a bounded
$q$ is a self-adjoint bounded operator acting in ${\cal H}^T$. The
last relation in (\ref{Delay rel}) can be written as $\Delta W^Tf
- W^T f_{tt}=q W^Tf$ that is just a form of writting the wave
equation (\ref{B1}). Taking into account the density of ${\cal
M}^T$ in ${\cal F}^T$, it is easy to conclude that a set of pairs
 \begin{equation}\label{Grapf q}
\left\{\langle \Delta W^Tf - W^T f_{tt}\,, W^T f\rangle\,|\,\,f
\in {\cal M}^T \right\}
 \end{equation}
determines the graph of the multiplication by $q$ and, hence,
determines the potential $q\big|_{\Omega^T}$.

\subsubsection*{Response operators}
$\bullet$\,\,\,In system $\alpha^T$, the input/output
correspondence is realized by a {\it response operator} $R^T: {\cal F}^T
\to {\cal F}^T$,
 $$
R^T f\,:=\,u^f\big|_{\Sigma^T}\,.
 $$
By the above-mentioned regularity of $u^f$, it acts continuously
from ${\cal F}^T$ to $H^{\frac{1}{5}-2\varepsilon}(\Sigma^T)$ and, hence,
is a compact operator. The following is some of its basic
properties. We use the auxiliary operators $Y^T, J^T: {\cal F}^T \to
{\cal F}^T$,
 $$
\left(Y^T f\right)(\cdot,t)\,:=\, f(\cdot, T-t)\,, \quad
\left(J^Tf\right)(\cdot,t)\,:=\,\int_0^t f(\cdot,s)\,ds\,, \qquad
0 \leqslant  t \leqslant T\,.
 $$
Note that $(Y^T)^*=(Y^T)^{-1}=Y^T$ and $(Y^T)^2={\mathbb I}_{{\cal F}^T}$
holds.

 \begin{Lemma}\label{L2}
For $T>0$ and $0 \leqslant \xi \leqslant T$, the relations
 \begin{equation}\label{R^T rel}
R^T {\cal T}^T_{T-\xi}\,=\,{\cal T}^T_{T-\xi}R^T\,; \quad R^T
J^T\,=\,J^T R^T\,; \quad (Y^T R^T)^*\,=\,Y^T R^T
 \end{equation}
are valid.
\end{Lemma}
\noindent$\square$\,\,\,\,The first relation follows from
(\ref{Delay rel}). The second is a simple consequence of the
first. Prove the third one.

Let controls $f,g$ belong to the smooth class ${\cal M}^T$, which
is dense in ${\cal F}^T$. Cauchy conditions (\ref{B2}) imply
 $$
u^f(\cdot,t)\big|_{t=0}=u^f_t(\cdot,t)\big|_{t=0}=u^g(\cdot,T-t)\big|_{t=T}=u^g_t(\cdot,T-t)\big|_{t=T}=0\,.
 $$
Also, since each $f \in {\cal M}^T$ vanishes near $t=0$, the wave
$u^f(\cdot,T)$ vanishes near $\Gamma^T$ by locality (\ref{supp
u-1}).

Integrating by parts, one has
 \begin{align*}
& 0=\int_{\Omega^T \times [0,T]}[u^f_{tt}-\Delta
u^f+qu^f](x,t)\,u^g(x,T-t)\,dx\,dt=\\
& =\int_{\Sigma^T}[u^f(\gamma,t)\,\partial_\nu u^g(\gamma,T-t)-
\partial_\nu
u^f(\gamma,t)\,u^g(\gamma,T-t)]\,d\Gamma\,dt+\\
& +\int_{\Omega^T \times [0,T]}u^f(x,t)[u^g_{tt}-\Delta
u^g+qu^g](x,T-t)\,dx\,dt=\\
&
\overset{(\ref{B3})}=\int_{\Sigma^T}[u^f(\gamma,t)\,g(\gamma,T-t)-
f(\gamma,t)\,u^g(\gamma,T-t)]\,d\Gamma\,dt=\\
& =(R^Tf,Y^Tg)_{{\cal F}^T}-(f,Y^T R^Tg)_{{\cal F}^T}=(Y^T R^Tf,g)_{{\cal F}^T}-(f,Y^T
R^Tg)_{{\cal F}^T}.
 \end{align*}
Thus, we have $(Y^T R^Tf,g)_{{\cal F}^T}=(f,Y^T R^Tg)_{{\cal F}^T}$. Since
${\cal M}^T$ is dense in ${\cal F}^T$, we get the last equality in
(\ref{R^T rel}).\,\,\,$\blacksquare$
\smallskip

\noindent$\bullet$\,\,\,There is one more object of system
$\alpha^T$ related with the input/output correspondence.

Denote $D^{2T}:={\rm in}\,\,\,\{(x,t)\,|\,\, x \in
{\Omega^T},\,\,t<2T-\tau(x)\}$. The problem
 \begin{align}
\label{C1}
 & u_{tt}-\Delta u+qu=0 &&{\rm in}\,\,\,D^{2T}\\
\label{C2}
 & u\big|_{t<\tau(x)}\,=\,0 &&{\rm in}\,\,\,\overline{D^{2T}}\\
\label{C3} &  \partial_\nu u=f \quad &&{\rm
on}\,\,\,\overline{\Sigma^{2T}}\,,
\end{align}
can be regarded as a natural extension of problem
(\ref{B1})--(\ref{B3}). Such an extension does exist and is well
posed  owing to the finiteness of the domains of influence
(hyperbolicity). Its solution $u^f$ is determined by
$q\big|_{\Omega^T}$.

With problem (\ref{C1})--(\ref{C3}) one associates an {\it
extended response operator} $R^{2T}: {\cal F}^{2T}\to {\cal
F}^{2T}$,
 $$
R^{2T} f\,:=\,u^f\big|_{\Sigma^{2T}}\,.
 $$
It is a compact operator with the properties quite analogous to
(\ref{R^T rel}):
 \begin{align}
\notag & R^{2T} {\cal T}^{2T}_{2T-\xi}\,=\,{\cal
T}^{2T}_{2T-\xi}R^{2T}\,, \quad 0 \leqslant \xi \leqslant 2T\,; \quad R^{2T}
J^{2T}\,=\, J^{2T} R^{2T}\,;\\
\label{R^2T rel} &  (Y^{2T} R^{2T})^*\,=\,Y^{2T} R^{2T}\,.
 \end{align}
Along with the solution $u^f$, operator $R^{2T}$ is determined by
$q\big|_{\Omega^T}$. By the latter, this operator must be regarded
as an intrinsic object of system $\alpha^T$ (but not
$\alpha^{2T}$). Note in addition that $R^{2T}$ is meaningful at a
very general level: see \cite{DSBC}.

\subsubsection*{Connecting operator}
$\bullet$\,\,\,A key object of the BC-method is a {\it connecting
operator} $C^T: {\cal F}^T \to {\cal F}^T$,
 \begin{equation}\label{def C^T}
C^T:=(W^T)^\ast W^T\,.
 \end{equation}
By the definition, we have
 \begin{equation*}
 (C^T f, g)_{{\cal F}^T} =
(W^T f, W^T g)_{{\cal H}^T} = \left(u^f(\cdot,T),
u^g(\cdot,T)\right)_{{\cal H}^T}\,,
 \end{equation*}
 i.e., $C^T$ connects the Hilbert metrics of the outer and
inner spaces. It is a compact (because $W^T$ is) and nonnegative
operator: $(C^Tf,f)_{{\cal F}^T}\geqslant 0$ holds for all $f \in {\cal F}^T$.
Moreover, since ${\rm Ker\,} C^T={\rm Ker\,} W^T$, Lemma \ref{L1} provides its
{\it positivity}:
 \begin{equation*}
 (C^T f, f)_{{\cal F}^T}\,>\,0 \qquad {\rm for}\,\,\,\,0\not= f \in {\cal F}^T\,,\,\,T<T_*.
 \end{equation*}
\smallskip

\noindent$\bullet$\,\,\,Recall that the image operator $ I^T$
introduced in section 1 acts from $L_2(\Omega^T)$ to
$L_2(\Sigma^T)$. In what follows we identify these spaces with
${\cal H}^T$ and ${\cal F}^T$ respectively, and regard $I^T$ as a
map from ${\cal H}^T$ to ${\cal F}^T$.

The definition of images easily implies $Y^T I^T{\cal H}^\xi
\subset {\cal F}^{T,\,\xi}$, whereas (\ref{W^T FTxis subset Hxis})
(for $\sigma=\Gamma$) provides $Y^T I^T W^T {\cal F}^{T,\,\xi}
\subset {\cal F}^{T,\,\xi}$. The latter means that an operator
$Y^T I^T W^T$ is {\it triangular} with respect to the family of
subspaces ({\it nest}) $\{{\cal F}^{T,\,\xi}\}_{0 \leqslant \xi
\leqslant T}$ \,\cite{Dav}.

For the connecting operator, the relations
 \begin{equation}\label{C^T triang factor}
C^T\,\overset{(\ref{def C^T})}=\,(W^T)^* W^T\,\overset{(\ref{I^T
properties})}=\,(Y^T I^TW^T)^*(Y^T I^T W^T)
 \end{equation}
hold and show that operator $Y^T I^T W^T$ provides  a {\it
triangular factorization} of the connecting operator with respect
to the nest $\{{\cal F}^{T,\,\xi}\}_{0 \leqslant \xi \leqslant T}$
\cite{GohKr, Dav}.

\smallskip

\noindent$\bullet$\,\,\, A significant fact is that the connecting
operator is determined by the extended response operator via an
explicit formula:
\begin{equation}\label{C^T via R^2T}
 C^T\,=\,-\,\frac{1}{2}~(S^T)^* R^{2T} J^{2T}
S^T\,,
\end{equation}
where the map $S^T: {\cal F}^T \to {\cal F}^{2T}$ extends the
controls from $\Sigma^T$ to $\Sigma^{2T}$ by oddness:
 \begin{equation*}
 \left(S^T f\right)(\cdot,t)\,=\,
          \begin{cases}
   f(\cdot,t)\,, & 0 \leqslant t < T\\
   -f(\cdot, 2T-t)\,, & T \leqslant t \leqslant 2T\,.
          \end{cases}
 \end{equation*}
In \cite{BIP97, BIP07}, a relevant analog of this representation
is proved for the case of the Dirichlet boundary controls. To
modify the proof for obtaining (\ref{C^T via R^2T}) needs just a
minor correction.

\subsection{System $\alpha^T_*$}\label{subsec System alpha^T*}
A dynamical system associated with the problem
\begin{align}
\label{D1}
 & v_{tt}-\Delta v+qv=0 &&{\rm in}\,\,\,\{(x,t)\,|\,\, x \in {\Omega^T},\,\,t>\tau(x)\}\\
\label{D2}
 & v\big|_{t=T}=0\,,\,\,\,v_t\big|_{t=T}= y \in {\cal H}^T\\
\label{D3}
 &  \partial_\nu v=0 \quad &&{\rm on}\,\,\,\Sigma^{T}
\end{align}
is denoted by $\alpha^T_*$ and said to be {\it dual} to system
$\alpha^T$. Its solution $v=v^y(x,t)$ describes a wave, which is
initiated by the velocity perturbation $y$ and propagates (in the
reversed time) in $\Omega$. The problem is well posed owing to the
finiteness of the domain of influence property.

Integration by parts provides the well-known relation
 \begin{equation*}
(u^f(\cdot,T), y)_{{\cal H}^T}\,=\,(f, v^y)_{{\cal F}^T}\,,\qquad f \in
{\cal F}^T,\,\,y\in {\cal H}^T\,.
 \end{equation*}
It is the relation, which motivates the term `dual'\,\cite{BIP97,
BIP07}.

In the dual system, the state/observation correspondence is
realized by an {\it observation operator} $O^T: {\cal H}^T \to {\cal F}^T$,
 $$
O^T y\,:=\,v^y \big |_{\Sigma^T}\,.
 $$
Being written in the form $(W^Tf, y)_{{\cal H}^T}\,=\,(f,
O^Ty)_{{\cal F}^T}$, the duality relation  leads to the equality
 \begin{equation}\label{O=W^*}
O^T\,=\,(W^T)^*\,.
 \end{equation}
It implies ${\rm Ker\,} O^T={\cal H}^T\ominus \overline{{\rm
Ran\,} W^T}$, whereas (\ref{W^T FTxis = Hxis}) (for
$\sigma=\Gamma$) follows to the equality ${\rm Ker\,} O^T=\{0\}$.
The latter is interpreted as a {\it boundary observability} of the
dual system.

\section{Visualization of waves}\label{sec Visualization}

\subsection{Devices}\label{subsec Devices}

\subsubsection*{Propagation of jumps in $\alpha^T_*$}
A very general fact of the propagation of singularities theory for
the hyperbolic equations is that discontinuous data produce
discontinuous solutions, the discontinuities propagating along
bicharacteristics and being supported on characteristic surfaces.
Here we deal with the Cauchy problem (\ref{D1})--(\ref{D3}) with a
$y$ having jumps of special kind. Our goal is to describe the
corresponding jumps of the image $O^T y$. The description is
provided by the proper Geometrical Optics formulae.  Since the
GO-technique is rather cumbersome, we have to restrict ourselves
to heuristic considerations and references to our papers
\cite{BKach94, BIP97}, where the rigorous analysis is developed.

We start with a simpler  case $T<T_c$: the simplification is that
the surfaces $\Gamma^\xi$ are smooth as $\xi \leqslant T$. A
characteristic function (indicator) of a set $A$ is denoted by
$\chi_A$:
 $$
\chi_A(p)\,:=\,\begin{cases}
                 1\,,  &  p \in A\\
                 0\,,  &  p \not \in A
               \end{cases}\,.
 $$

\noindent$\bullet$\,\,\,\,Fix a $\xi$ and (small) $\Delta \xi$
provided $0<\xi<\xi+\Delta \xi <T$. A subdomain
 $$
\Delta \Omega^\xi\,:=\,\overline{\Omega^{\xi+\Delta \xi} \setminus
\Omega^\xi}\,\,\subset\, \Omega^T
 $$
is a thin layer between the smooth surfaces $\Gamma^{\xi+\Delta
\xi}$ and $\Gamma^\xi$.

Take a $y \in C^\infty(\overline{\Omega^T})$. A `slice'
$\chi_{\Delta \Omega^\xi}y$ is a piece-wise smooth function
supported in $\overline{\Delta \Omega^\xi}$. Generically, it has
the jumps at $\Gamma^\xi$ and $\Gamma^{\xi+\Delta \xi}$. In what
follows, the jump at $\Gamma^\xi$ is of our main interest, whereas
the jump at $\Gamma^{\xi+\Delta \xi}$ is introduced just for
technical convenience.

Return to system (\ref{D1})--(\ref{D3}). Putting
$v_t\big|_{t=T}=\chi_{\Delta \Omega^\xi}y$ in (\ref{D2}), we get a
Cauchy problem with discontinuous data. In particular, the data
have a jump at $\Gamma^\xi$:
 \begin{equation}\label{jump data}
v_t\left(x(\gamma,\tau),T\right)\bigg|_{\tau=\xi-0}^{\tau=\xi+0}=
\,y(x(\gamma,\xi))\,-\,0\,=\,y(x(\gamma,\xi))\,.
 \end{equation}
As a consequence, the solution $v^{\chi_{\Delta \Omega^\xi}y}$
turns out to be non-smooth. The following is some details specific
for problem (\ref{D1})--(\ref{D3}).
\smallskip

\begin{figure}\label{Fig}
\centering
\includegraphics[width=3.0in,height=2.5in]{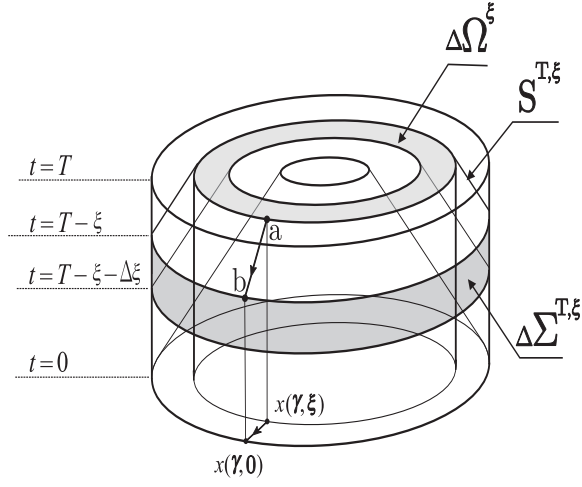}
\caption{Propagation of jump}
 \end{figure}

\noindent$\bullet$\,\,\,A velocity perturbation ${\chi_{\Delta
\Omega^\xi}y}$, which initiates the wave process, is separated
from the boundary with the distance $\xi$. Therefore, by the
finiteness of domain of influence principle, the solution
$v^{\chi_{\Delta \Omega^\xi}y}$ vanishes for $t>T-\xi - \tau(x)$,
i.e., over a characteristic surface $S^{T,\,\xi}:=\{(x,t)\in
\overline {\Omega^T} \times [0,T]\}$ (see \ref{Fig}).

\noindent$\bullet$\,\,\,Jumps of $v_t(\cdot,T)$ initiate jumps of
the velocity $v^{\chi_{\Delta \Omega^\xi}y}_t$. One of the
velocity jumps is located at the characteristic $S^{T,\,\xi}$
\footnote{another jumps also do occur but are beyond our
interest}. This jump propagates along the space-time rays
$r_\gamma^{T,\,\xi}$, which constitute the characteristic:
 \begin{align*}
 &r_\gamma^{T,\,\xi}:=\{(x,t)\in
\overline{\Omega^T}\times[0,T]\,|\,\,x=x(\gamma,\xi-\tau),\,\,t=T-\tau:\,\,0 \leqslant \xi \leqslant
T\}\,,\\
 & S^{T,\,\xi}\,=\,\bigcup_{\gamma \in \Gamma}
r_\gamma^{T,\,\xi}\,.
 \end{align*}
{\blue The} jump, which moves along $r_\gamma^{T,\,\xi}$, starts
from the point $a=(x(\gamma,\xi),T)$ and reaches the boundary at
$b=(x(\gamma,0),T-\xi)$. By (\ref{jump data}), at the `input' $a$
the value ({\it amplitude}) of the jump is $y(x(\gamma,\xi))$. At
the endpoint $b$, its amplitude is found by the GO-technique,
which provides
 \begin{equation}\label{jump v_t}
v^{\chi_{\Delta \Omega^\xi}y}_t(\left(
x(\gamma,0),t\right)\bigg|_{t=T-\xi-0}^{t=T-\xi+0}=0-
\beta^{\frac{1}{2}}(\gamma,\xi)y(x(\gamma,\xi))=-\,\beta^{\frac{1}{2}}(\gamma,\xi)y(x(\gamma,\xi))\,.
 \end{equation}
This relation corresponds to the well-known GO-law: the ratio of
the input and output jump amplitudes is governed by the factor
$\beta$, which is determined by the spreading of rays
$r^{T,\,\xi}_\gamma$ \cite{Ikawa, BKach94, BIP97}.
\smallskip

\noindent$\bullet$\,\,\,By the aforesaid, a trace $v^{\chi_{\Delta
\Omega^\xi}y}_t\big|_{\Sigma^T}$ vanishes on $\Gamma \times
(T-\xi,T]$ and has a jump at the cross-section $\Sigma^T\cap
S^{T,\,\xi} = \Gamma \times \{t=T-\xi\}$. In the mean time, by the
regularity results, this trace is continuous as an
$H^{\frac{1}{2}}(\Gamma)$-valued function of $t \in [0,T-\xi]$
\footnote{this property can be derived from Theorem 3.3 of
\cite{LTr-2}.}. The following considerations specify the behavior
of $v^{\chi_{\Delta \Omega^\xi}y}_t\big|_{\Sigma^T}$ near (and
below) this cross-section.

Let
 $$
\Delta \Sigma^{T,\,\xi}:=\{(\gamma,t)\in \Sigma^T\,|\,\,\gamma \in
\Gamma,\,\,T-\xi-\Delta\xi \leqslant t \leqslant T-\xi\}
 $$
be a thin `belt' near the cross-section (see Fig. 1),
$\chi_{\Delta \Sigma^{T,\,\xi}}$ its indicator. A function on
$\Sigma^T$ of the form $\chi_{\Delta
\Sigma^{T,\,\xi}}\left[v^{\chi_{\Delta
\Omega^\xi}y}_t\big|_{\Sigma^T}\right]$ is a `slice' of the
boundary trace of the velocity. By (\ref{jump v_t}), one can
represented it as
 \begin{align}
\notag & \left(\chi_{\Delta \Sigma^{T,\,\xi}}\left[v^{\chi_{\Delta
\Omega^\xi}y}_t\big|_{\Sigma^T}\right]\right)(\gamma,t)\,=\\
\label{slice of velocity} &        =
       \begin{cases}
-\beta^{\frac{1}{2}}(\gamma,\xi)y(x(\gamma,\xi))+w^{\xi,\Delta\xi}(\gamma,t)\,,   & (\gamma,t)\in \Delta \Sigma^{T,\,\xi} \\
0\,,          &(\gamma,t)\in \Sigma^T\setminus\Delta
\Sigma^{T,\,\xi}
        \end{cases}\,\,,
 \end{align}
where the first summand in the first line {\it does not depend on
$t$} and, hence, obeys $\|\beta^{\frac{1}{2}}y\|^2_{L_2(\Delta
\Sigma^{T,\,\xi})}\thicksim \Delta \xi$, whereas the second
summand satisfies\\ $\|w^{\xi,\Delta\xi}\|^2_{L_2(\Delta
\Sigma^{T,\,\xi})}\thicksim o(\Delta \xi)$ uniformly with respect
to $\xi \in [0,T]$ and (small enough) $\Delta \xi>0$
\cite{BKach94, BIP97}. So, the first summand is dominating.

\subsubsection*{Amplitude integral}
$\bullet$\,\,\,Choose a partition
$\Xi=\{\xi_i\}_{i=0}^N:\,\,\,0=\xi_0<\xi_1< \dots<\xi_N=T$ of the
segment $[0,T]$ and denote
 \begin{align*}
& \Delta \xi_i=\xi_i-\xi_{i-1}\,, \quad
\Delta\Sigma^{T,\,\xi_i}=\Gamma \times [T-\xi_i-\Delta\xi_i\,,
T-\xi_i]\,, \quad \Delta
\Omega^{\xi_i}=\overline{\Omega^{\xi_i}\setminus\Omega^{\xi_{i-1}}}\,,\\
& i=1,2,\dots N\,\,\,\,\,\,\,(\Omega^0:=\emptyset); \qquad
r_\Xi\,=\,\max_{i=1,\dots,N}\Delta \xi_i\,.
 \end{align*}
Summing up the terms of the form (\ref{slice of velocity}) and
recalling the definition of images, we get
 \begin{align}
\notag &\left(\sum\limits_{i=1}^N\chi_{\Delta
\Sigma^{T,\,\xi_i}}\left[v^{\chi_{\Delta
\Omega^{\xi_i}}y}_t\big|_{\Sigma^T}\right]\right)(\gamma,T-t)\,=\\
& =\,-\,\left(I^Ty\right)(\gamma,t) +\delta^{y,\,\Xi}(\gamma,t),
\qquad (\gamma,t)\in \Sigma^T\,,\label{integral sums}
 \end{align}
where $\|\delta^{y,\,\Xi}\|_{L_2(\Sigma^{T})}\to 0$ as $r_\Xi \to
0$. Substituting $t$ by $T-t$, we see that, for the given smooth
$y \in {\cal H}^T$, the sums converge to $-Y^T I^T y$ by the norm
in ${\cal F}^T$. The smallness of $\delta^{y,\,\Xi}$ is justified by
perfect analogy with the case of the problem with Dirichlet
boundary controls \cite{BKach94, BIP97}.
\smallskip

\noindent$\bullet$\,\,\,Here we interpret (\ref{integral sums}) in
operator terms.

Let $X^{T,\,\xi}$ be a projection in ${\cal F}^T$ onto ${\cal
F}^{T,\,\xi}$, which cuts off controls onto $\Gamma \times
[T-\xi,T]$. The difference $\Delta X^{T,\,\xi_i}=
X^{T,\,\xi_i}-X^{T,\,\xi_{i-1}}$ is also the projection cutting
off controls onto the belt $\Delta \Sigma^{\xi_i,\,T}$: $\Delta
X^{T,\,\xi_i}f=\chi_{\Delta \Sigma^{T,\,\xi_i}}f$.

By $G^\xi$ we denote a projection in ${\cal H}^T$ onto ${\cal
H}^\xi$, which cuts off functions onto $\Omega^\xi$. The
difference $\Delta G^{\xi_i}= G^{\xi_i}-G^{\xi_{i-1}}$ cuts off
functions onto the layer $\Delta \Omega^{\xi_i}$: $\Delta
G^{\xi_i}y=\chi_{\Delta \Omega^{\xi_i}}y$.

Recalling the definition of the observation operator, one can
represent the summands in (\ref{integral sums}) as
 $$
\chi_{\Delta \Sigma^{T,\,\xi_i}}\left[v^{\chi_{\Delta
\Omega^{\xi_i}}y}_t\big|_{\Sigma^T}\right]=\Delta
X^{T,\,\xi_i}\partial_tO^T\Delta G^{\xi_i}y
 $$
and then write (\ref{integral sums}) in the form
 \begin{equation}\label{AI-1}
\lim \limits_{r_\Xi \to 0}\left[\sum\limits_{i=1}^N\Delta
X^{T,\,\xi_i}\partial_tO^T\Delta
G^{\xi_i}\right]y=:\left[\int_{[0,T]}d X^{T,\,\xi}\,\partial_t
O^T\,dG^{\xi}\right]y = Y^T I^T y\,.
 \end{equation}
An operator construction in the square brackets is said to be an
{\it amplitude integral} (AI). It represents the image of $y$ as a
collection of the wave jumps, which pass through $\Omega^T$ and
are detected by the external observer at the boundary.
\smallskip

\noindent$\bullet$\,\,\,Recall that (\ref{AI-1}) is derived under
the assumption $T<T_c$. The case $T>T_c$ is more complicated since
the equidistant surfaces $\Gamma^\xi$ can be non-smooth and
disconnected. However, a remarkable fact is that representation
(\ref{AI-1}) does survive: it is valid for {\it any} $T<T_*$. For
the system $\alpha^T$ with Dirichlet boundary controls, this
result is stated in \cite{BKach94, BIP97}. To modify it for the
case of Neumann controls requires just a minor technical changes.
So, the following does occur.
 \begin{Proposition}\label{Prop 1}
For any positive $T<T_*$, the sums in (\ref{AI-1}) converge to the
limit
 \begin{equation}\label{AI-2}
\lim \limits_{r_\Xi \to 0}\sum\limits_{i=1}^N\Delta
X^{T,\,\xi_i}\partial_tO^T\Delta G^{\xi_i}=:\int_{[0,T]}d
X^{T,\,\xi}\,\partial_t O^T\,dG^{\xi} = Y^T I^T
 \end{equation}
in the weak operator topology.
 \end{Proposition}

\subsubsection*{$W^T$ via amplitude integral}
\noindent$\bullet$\,\,\,Multiplying (\ref{AI-2}) by $W^T$ from the
right, we get an operator $V^T: {\cal F}^T \to {\cal F}^T$,
 \begin{equation}\label{V^T def}
V^T:=Y^T I^T W^T=\left[\int_{[0,T]}d X^{T,\,\xi}\,\partial_t
O^T\,dG^{\xi}\right]W^T,
 \end{equation}
which satisfies
\begin{equation}\label{V^T factor}
V^T {\cal F}^{T,\,\xi}\,\subset\,{\cal F}^{T,\,\xi}\,, \qquad
(V^T)^*V^T\overset{(\ref{C^T triang factor})}=C^T\,.
 \end{equation}
Thus, $V^T$ provides triangular factorization of the connecting
operator with respect to the nest $\{{\cal F}^{T,\,\xi}\}_{0
\leqslant \xi \leqslant T}$.

\smallskip

\noindent$\bullet$\,\,\, Any densely defined closable linear
operator acting from a Hilbert space to a Hilbert space can be
represented in the form of a {\it polar decomposition} (see, e.g.,
\cite{BirSol}). For the control operator, such a decomposition is
 \begin{equation}\label{W=U|W|}
W^T\,=\,U^T
|W^T|\,:=\,U^T\left[(W^T)^*W^T\right]^{\frac{1}{2}}\,\overset{(\ref{C^T
triang factor})}=\,U^T\left[C^T\right]^{\frac{1}{2}}\,,
 \end{equation}
where $|W^T|: {\cal F}^T \to {\cal F}^T$ is a modulo of $W^T$, and $U^T: {\cal F}^T \to
{\cal H}^T$ is an isometry, which maps ${\rm Ran\,} |W^T| \subset {\cal F}^T$ onto
${\rm Ran\,} W^T \subset {\cal H}^T$ by the rule
 \begin{equation}\label{U^T|W^T|f=W^T f}
U^T|W^T|f \,=\, W^T f\,, \qquad f \in {\cal F}^T\,.
 \end{equation}
By (\ref{W^T FTxis = Hxis}) with $\sigma=\Gamma$, for {\it any}
$T>0$ one has $\overline{{\rm Ran\,} W^T}= {\cal H}^T$. In the mean time, for
$T<T_*$, we have
 $$
\overline{{\rm Ran\,} |W^T|}={\cal F}^T\ominus {\rm Ker\,} |W^T|={\cal F}^T\ominus {\rm Ker\,}
W^T\overset{{\rm Lemma}\, \ref{L1} }={\cal F}^T\,.
 $$
As a result, if $T<T_*$ then $U^T$ can be extended by continuity
from ${\rm Ran\,} |W^T|$ to ${\cal F}^T$, the extension being a
{\it unitary} operator, which maps ${\cal F}^T$ onto ${\cal H}^T$.
In what follows, we assume that such an extension is done; it
satisfies
 \begin{equation}\label{U^*U=1}
(U^T)^* U^T\,=\,{\mathbb I}_{{\cal F}^T}\,, \qquad
U^T(U^T)^*\,=\,{\mathbb I}_{{\cal H}^T}\,.
 \end{equation}
\smallskip

\noindent$\bullet$\,\,\,Recall that $G^\xi$ projects in ${\cal
H}^T$ onto ${\cal H}^\xi$. We say a projection $P^\xi$ in ${\cal
H}^T$ onto the subspace $\overline{W^T{\cal F}^{T,\,\xi}}$ (formed
by waves) to be a {\it wave projection}. A crucial point of our
approach is the equality
 \begin{equation}\label{P^xi=G^xi}
P^\xi\,\overset{(\ref{W^T FTxis = Hxis})}=\, G^\xi\,, \qquad 0 \leqslant
\xi \leqslant T\,,
 \end{equation}
which corresponds to the controllability of system $\alpha^T$.

Let $\tilde P^{T,\,\xi}$ be  a projection in ${\cal F}^T$ onto the
subspace $\overline{|W^T|{\cal F}^{T,\,\xi}}$. By (\ref{U^T|W^T|f=W^T f}), one
has
 \begin{equation}\label{U^T tilde P^xi=P^xi}
U^T \tilde P^{T,\,\xi}\,=\,P^\xi U^T\,, \qquad 0 \leqslant \xi \leqslant T
 \end{equation}
that implies
 \begin{align}
\notag & O^T G^\xi W^T\overset{(\ref{O=W^*}), (\ref{P^xi=G^xi})} =
(W^T)^* P^\xi W^T \overset{(\ref{W=U|W|})} =|W^T|(U^T)^*P^\xi
U^T|W^T|=\\
\label{Auxilliary} & \overset{(\ref{U^T tilde
P^xi=P^xi})}=|W^T|\,\tilde P^{T,\,\xi}\,|W^T|
 \end{align}
for $0 \leqslant \xi \leqslant T$.
\smallskip

\noindent$\bullet$\,\,\,Multiplying equality (\ref{V^T def}) by
the isometry $(I^T)^*Y^T$ from the left, and taking into account
(\ref{Auxilliary}), we get
 \begin{equation}\label{W^T=U^T|W^T|}
W^T=U^T|W^T|,\quad  U^T\,=\,(I^T)^*Y^T\left[\int_{[0,T]}d
X^{T,\,\xi}\,\partial_t |W^T|\,d\tilde P^{T,\,\xi}\right]\,.
 \end{equation}
Here the operators $I^T,\, Y^T, \,X^{T,\,\xi}$ are standard (do
not depend on potential $q$), whereas projections $\tilde
P^{T,\,\xi}$ are obviously determined by $|W^T|$. Operator $W^T$
is {\it triangular} with respect to the pair of the nests $\{{\cal
F}^{T,\,\xi}\}$ and $\{{\cal H}^\xi\}$ that means $W^T {\cal
F}^{T,\,\xi} \subset {\cal H}^\xi,\,\,\,0 \leqslant \xi \leqslant
T$ (see (\ref{W^T FTxis = Hxis})). From the operator theory
viewpoint, representation (\ref{W^T=U^T|W^T|}) enables one to
recover a triangular operator $W^T$ via its modulo $|W^T|$, the
`phase' part $U^T$ being expressed via a relevant operator
integral. The integral into the square brackets is referred to as
a {\it diagonal} of operator $\partial_tW^T$ with respect to the
nests $\{{\cal F}^{T,\,\xi}\}$ and $\{{\cal H}^\xi\}$ \cite{Dav,
BPush}.
\smallskip

\noindent$\bullet$\,\,\,Introduce an operator $A^T: {\cal F}^T \to {\cal F}^T$
by
 \begin{equation}\label{A^T def}
A^T\,:=\,Y^T\int_{[0,T]}d X^{T,\,\xi}\,\partial_t
[C^T]^{\frac{1}{2}}\,d\tilde P^{T,\,\xi}\,.
 \end{equation}
With regard to (\ref{P^xi=G^xi}) and (\ref{U^T tilde P^xi=P^xi}),
one can write (\ref{AI-2}) in the form $A^T (U^T)^*=I^T$ that
enables one to represent the phase operator in the form
 \begin{equation*}
U^T\overset{(\ref{W^T=U^T|W^T|})}=(I^T)^* A^T\,.
 \end{equation*}
By (\ref{U^*U=1}) and (\ref{I^T properties}), this representation
implies
 \begin{equation}\label{U=I^*YA}
 \quad
(A^T)^*A^T\,=\,{\mathbb I}_{{\cal F}^T},\qquad
A^T(A^T)^*\,=\,G_{\Theta^T}\,.
 \end{equation}
Now, writing (\ref{W^T=U^T|W^T|}) in the form
 \begin{equation}\label{Basic Repres}
W^T\,=\,(I^T)^* A^T\,[C^T]^{\frac{1}{2}}\,,
 \end{equation}
we obtain the representation of the control operator, which plays
a basic role in solving inverse problems. The reason is the
following.

Operator $R^{2T}$ formalizes information, which the external
observer gets from measurements at the boundary $\Gamma$. The
waves $u^f$ propagate into $\Omega$ and are {\it invisible} for
him. However, the observer can determine $C^T$ via (\ref{C^T via
R^2T}), find $[C^T]^{\frac{1}{2}}$, construct the integral
(\ref{A^T def}), determine $W^T$ via (\ref{Basic Repres}), and
eventually recover invisible waves $u^f(\cdot,T)=W^Tf$. In the
BC-method, such a remarkable option is referred to as a {\it
visualization of waves}.

\subsection{Solving the inverse problem}

\subsubsection*{Setup}
As is mentioned in section \ref{subsec System alpha^T}, the
extended response operator $R^{2T}$ depends on the potential
locally: it is determined by $q\big|_{\Omega^T}$. Such a locality
motivates the following setup of the inverse problem.
\smallskip

\noindent{\bf (IP)}\,\,\,{\it Given operator $R^{2T}$, to recover
 potential $q$ in the subdomain $\Omega^T$}.
\smallskip

\noindent The IP will be solved for an arbitrary fixed  $T<T_*$.
Surely, such an option enables one to determine $q$ in the whole
$\Omega$ if $R^{2T}$ is given for a $T \geqslant T_*$.

\subsubsection*{Procedure}
Preparatory to solving the IP, recall that geometry of the wave
propagation in system $\alpha^T$ is governed by the leading part
$\partial^2_t - \Delta$ of the wave equation (\ref{A1}). Since
this part does not depend on the potential, the geometry is {\it
Euclidean}.  Therefore, we have the right to regard all the
geometric objects and parameters ($\Omega^\xi$, sgc, $\Theta^T$,
$\beta$, $T_*$, etc) as {\it known} and use them for determination
of $q$. In particular, we can use the image operator $I^T$.
\smallskip

Let $T<T_*$ be fixed. Given $R^{2T}$ one can recover $q$ in
$\Omega^T$ by the following procedure.
\smallskip

\noindent{\it Step 1.}\,\,\,Find $C^T$ by (\ref{C^T via R^2T}).
Determine $[C^T]^{\frac{1}{2}}$.
\smallskip

\noindent{\it Step 2.}\,\,\,Determine the subspaces
$[C^T]^{\frac{1}{2}}{\cal F}^{T,\,\xi}$ and the corresponding projections
$\tilde P^{T,\,\xi}$  for $0 \leqslant \xi \leqslant T$.
\smallskip

\noindent{\it Step 3.}\,\,\,Construct the integral (\ref{A^T def})
and, then, recover $W^T$ via (\ref{Basic Repres}).
\smallskip

\noindent{\it Step 4.}\,\,\,Determine $q\big|_{\Omega^T}$ from the
graph (\ref{Grapf q}).
\smallskip

The IP is solved.

\subsection{Characterization of data}
\subsubsection*{Main result}
In addition to the procedure, which solves the IP, we provide the
necessary and sufficient conditions for its solvability.
 \begin{Theorem}\label{T1}
Let $0<T<T_*$. An operator ${\cal R}^{2T}: {\cal F}^{2T}\to {\cal
F}^{2T}$ is the extended response operator of {\it a} system
$\alpha^T$ with potential of the class $L_\infty(\Omega^T)$ if and
only if it satisfies the following conditions:
 \begin{enumerate}
\item ${\cal R}^{2T}$ is a compact operator obeying
 \begin{equation}\label{cal R}
Y^{2T}{\cal R}^{2T}=({\cal R}^{2T}Y^{2T})^*; \,\,\, {\cal
R}^{2T}{\cal T}^{2T}_{2T-\xi} = {\cal T}^{2T}_{2T-\xi}{\cal
R}^{2T}, \quad 0 \leqslant \xi \leqslant 2T\,.
 \end{equation}

\item An operator ${\cal C}^T: {\cal F}^T \to {\cal F}^T$,
 \begin{equation}\label{cal C^T}
{\cal C}^T\,:=\,-\,\frac{1}{2}\,(S^T)^*{\cal R}^{2T}J^{2T} S^T
 \end{equation}
is symmetric and positive: $({\cal C}^T f,f)_{{\cal F}^T}>0$ for $0 \not=
f \in {\cal F}^T$.

\item Let $\tilde{\cal P}^{T,\,\xi}$ be a projection in ${\cal
F}^T$ onto $\overline{[{\cal C}^T]^{\frac{1}{2}}{\cal
F}^{T,\,\xi}}$. An operator integral ${\cal A}^T: {\cal F}^T \to
{\cal F}^T$,
 \begin{equation}\label{cal U^T}
{\cal A}^T\,:=\,Y^T \int_{[0,T]}d X^{T,\,\xi}\,\partial_t [{\cal
C}^T]^{\frac{1}{2}}\,d\tilde {\cal P}^{T,\,\xi}
 \end{equation}
converges in the weak operator topology to an isometry, which
satisfies
 \begin{equation}\label{A^T* A^T=1}
({\cal A}^T)^*  {\cal A}^T\,=\,{\mathbb I}_{{\cal F}^T}, \quad {\cal A}^T
({\cal A}^T)^*\,=\,G_{\Theta^T}\,.
 \end{equation}
 \item An operator ${\cal W}^T: {\cal F}^T \to {\cal H}^T$
 \begin{equation}\label{Basic for cal W^T}
{\cal W}^T\,:=\,(I^T)^*{\cal A}^T [{\cal C}^T]^{\frac{1}{2}}
 \end{equation}
satisfies ${\cal W}^T{\cal M}^T \subset H^2(\Omega^T)$.

 \item
The relation
  \begin{equation}\label{boundary condition}
\partial_\nu{\cal W}^T f \big|_{\Gamma}\,=\,f (\cdot, T)\,, \qquad
f \in {\cal M}^T
 \end{equation}
is valid.

 \item The relation
 \begin{equation}\label{cal W^T new}
\overline{{\cal W}^T {\cal F}^{T,\,\xi}_\sigma}\,=\, {\cal H}^\xi_\sigma\,,\qquad 0 \leqslant \xi \leqslant T
 \end{equation}
holds for any open $\sigma \subseteq \Gamma$.
 \item The relation
\begin{equation}\label{Estimate Q^T}
\sup \limits_{0 \not=f\in {\cal M}^T}\frac{\| \Delta {\cal W}^T
f-{\cal W}^T f_{tt} \|_{{\cal H}^T}}{\|{\cal W}^T f\|_{{\cal H}^T}}\,<\,\infty
 \end{equation}
holds.

 \end{enumerate}
 \end{Theorem}

The proof consists of two parts.

\subsubsection*{Part I (necessity)}
$\square$\,\,\,Let ${\cal R}^{2T}=R^{2T}$, where $R^{2T}$ is the
extended response operator of a system $\alpha^T$ with potential
$q\in L_\infty(\Omega^T)$. The system possesses the connecting,
control, and phase operators $C^T$, $W^T$, and $U^T$ respectively.
\smallskip

\noindent {\it 1.}\,\,\,Relations (\ref{cal R}) hold by (\ref{R^2T
rel}).
\smallskip

\noindent $\it 2.$\,\,\, In view of (\ref{C^T via R^2T}), operator
${\cal C}^T$ defined by (\ref{cal C^T}) coincides with $C^T$,
which is a compact positive operator.
\smallskip

\noindent $\it 3.$\,\,\,The equality ${\cal C}^T=C^T$ implies
$\tilde {\cal P}^{T,\,\xi}=\tilde P^{T,\,\xi}$. Comparing
(\ref{A^T def}) with (\ref{cal U^T}), we conclude that ${\cal
A}^T=A^T$. Hence, (\ref{A^T* A^T=1}) follows from (\ref{U=I^*YA}).
\smallskip

\noindent $\it 4.$\,\,\, Comparing (\ref{Basic for cal W^T}) with
(\ref{Basic Repres}), we see that ${\cal W}^T$ coincides with
$W^T$. Hence, ${\cal W}^T{\cal M}^T \subset H^2(\Omega^T)$ holds
by the regularity results on the problem (\ref{A1})--(\ref{A3})
(see section \ref{subsec IBV-problem}).
\smallskip

\noindent $\it 5.$\,\,\,Since ${\cal W}^T=W^T$, the equality
(\ref{boundary condition}) is just a form of writing (\ref{B3}).
\smallskip

\noindent $\it 6.$\,\,\,(\ref{cal W^T new}) holds by (\ref{W^T
FTxis = Hxis}).
\smallskip

\noindent $\it 7.$\,\,\,Since ${\cal W}^T f=W^T f = u^f(\cdot,T)$,
we have
 \begin{align*}
& -\Delta {\cal W}^T f+{\cal W}^T f_{tt}=-\Delta
u^f(\cdot,T)+u^{f_{tt}}(\cdot,T)\overset{(\ref{Delay
rel})}=\\
& =- \Delta
u^f(\cdot,T)+u^f_{tt}(\cdot,T)\,\overset{(\ref{B1})}=\,q
u^f(\cdot,T)\,.
 \end{align*}
The inequality (\ref{Estimate Q^T}) is a consequence of $q \in
L_\infty(\Omega^T)$. \,\,\,$\blacksquare$

\subsubsection*{Part II (sufficiency)}
The proof of sufficiency is constructive: given ${\cal R}^{2T}$ we
provide a system $\alpha^T$ with the response operator
$R^{2T}={\cal R}^{2T}$. In fact, the construction follows the
procedure {\it Step 1-4}, which solves the IP.
\smallskip

$\square$\,\,\,Assume that ${\cal R}^{2T}$ obeys $\it 1$-$\it 5$.
\smallskip

\noindent$\bullet$\,\,\, Determine operator ${\cal C}^T$ by
(\ref{cal C^T}) and find $[{\cal C}^T]^{\frac{1}{2}}$. The latter
is also positive and injective.

Construct the operator integral in (\ref{cal U^T}) and get
operator ${\cal A}^T$. By  (\ref{A^T* A^T=1}), ${\cal A}^T$ is an
isometry in ${\cal F}^T$ with the range $G_{\Theta^T}{\cal F}^T$.
Hence, it satisfies $G_{\Theta^T}{\cal A}^T={\cal A}^T$.

Introduce operator ${\cal W}^T: {\cal F}^T \to {\cal H}^T$ in
accordance with (\ref{Basic for cal W^T}). Obviously, it is
injective. By (\ref{cal W^T new}) (for $\xi=T$ and
$\sigma=\Gamma$), its range ${\cal W}^T{\cal F}^T$ is dense in
${\cal H}^T$. Also, it satisfies
 \begin{align}
\notag & ({\cal W}^T)^*{\cal W}^T = [{\cal
C}^T]^{\frac{1}{2}}({\cal A}^T)^*I^T(I^T)^*{\cal A}^T[{\cal
C}^T]^{\frac{1}{2}}\overset{(\ref{I^T properties})}=[{\cal
C}^T]^{\frac{1}{2}}({\cal A}^T)^*G_{\Theta^T}{\cal A}^T[{\cal
C}^T]^{\frac{1}{2}}=\\
\label{{cal W}^*{cal W}= cal C} & = [{\cal
C}^T]^{\frac{1}{2}}({\cal A}^T)^*{\cal A}^T[{\cal
C}^T]^{\frac{1}{2}}\,\overset{(\ref{A^T* A^T=1})}=\,{\cal C}^T\,.
 \end{align}

\smallskip

\noindent$\bullet$\,\,\,Since ${\cal W}^T$ is injective, the set
of pairs
 $
\left\{\langle{\cal W}^Tf,\,{\cal W}^Tf_{tt}\rangle\,|\,\,f \in
{\cal M}^T\right\}
 $
constitutes the graph of a linear operator acting in ${\cal H}^T$. This
operator is denoted by $L^T: {\cal W}^Tf \mapsto {\cal
W}^Tf_{tt}$. It acts in ${\cal H}^T$ and is densely defined (on ${\cal
W}^T{\cal F}^T$).

Recall that the class of smooth controls ${\cal M}^T$ is dense in
${\cal F}^T$, its elements vanishing near $t=0$. The subclass
 $$
{\cal M}^T_0:=\{f \in {\cal M}^T\,|\,\,f \,\,{\rm
vanishes\,\,near\,\,} t=T\}
 $$
is also dense in ${\cal F}^T$. Hence, ${\cal W}^T{\cal M}^T_0$ is dense
in ${\cal H}^T$ by (\ref{cal W^T new}) for $\sigma=\Gamma,\, \xi=T$. As a
result, an operator $L^T_0:=L^T\big|_{{\cal W}^T{\cal M}^T_0}$ is
densely defined in ${\cal H}^T$. Show that it is symmetric.

Take $f, g \in {\cal M}^T_0$. Note that $S^Tf$ and $S^T g$ are
twice differentiable with respect to $t$ and vanish near $t=0$ and
$t=2T$. Also, note that the second relation in (\ref{cal R})
implies the commutation ${\cal
R}^{2T}\partial^2_t=\partial^2_t{\cal R}^{2T}$. Then, we have
 \begin{align*}
 & (L^T_0 {\cal W}^T f, {\cal W}^T g)_{{\cal H}^T}= (L^T{\cal W}^T f,
{\cal W}^T g)_{{\cal H}^T}= ({\cal W}^T f_{tt}, {\cal W}^T
g)_{{\cal H}^T}\overset{(\ref{{cal W}^*{cal W}= cal C})}=\\
 &=({\cal C}^Tf_{tt}, g)_{{\cal F}^T}
\overset{(\ref{cal C^T})}=-\,\frac{1}{2}\,([{\cal R}^{2T}J^{2T}
S^T]f_{tt}, S^Tg)_{{\cal F}^{2T}}=\\
 & =-\,\frac{1}{2}\,([{\cal
R}^{2T}J^{2T} S^Tf]_{tt}, S^Tg)_{{\cal F}^{2T}}
\overset{\star}=-\,\frac{1}{2}\,({\cal R}^{2T}J^{2T}
S^Tf, [S^Tg]_{tt})_{{\cal F}^{2T}}=\\
 &=-\,\frac{1}{2}\,({\cal
R}^{2T}J^{2T} S^Tf, S^T g_{tt})_{{\cal F}^{2T}} =
-\,\frac{1}{2}\,((S^T)^*{\cal R}^{2T}J^{2T} S^Tf, g_{tt})_{{\cal
F}^{T}}=\\
 &=({\cal C}^Tf, g_{tt})_{{\cal
F}^T}\overset{(\ref{{cal W}^*{cal W}= cal C})} =({\cal W}^T f,
{\cal W}^T g_{tt})_{{\cal H}^T}= ({\cal W}^T f, L^T{\cal W}^T g)_{{\cal H}^T}=\\
 & =({\cal W}^T f, L^T_0 {\cal W}^T g)_{{\cal H}^T}.
 \end{align*}
In $(\star)$ we integrate by part with respect to time in
${\cal F}^T=L_2(\Sigma^T)$. So, $L^T_0$ is symmetric.
\smallskip

\noindent$\bullet$\,\,\, Owing to (\ref{Estimate Q^T}), operator
$Q^T:=\Delta -L^T$ defined on the dense set ${\cal W}^T{\cal
F}^T\subset {\cal H}^T$, is bounded. By this, we assume that $Q^T$
is extended to ${\cal H}^T$ by continuity.

Operator $Q^T$ is self-adjoint. Indeed, in view of (\ref{boundary
condition}), for $f \in {\cal M}^T_0$ one has $\partial_\nu{\cal
W}^T f\big|_{\Gamma}=f\big|_{t=T}=0$, i.e., elements of ${\cal
W}^T{\cal M}^T_0$ satisfy the homogeneous Neumann boundary
condition on $\Gamma$. By the latter, the Laplacian $\Delta$ is
symmetric on ${\cal W}^T{\cal M}^T_0$. Hence, $Q^T\big|_{{\cal
W}^T{\cal M}^T_0}=\Delta\big|_{{\cal W}^T{\cal M}^T_0}-L^T_0$ is
symmetric on a dense set. Since it is bounded, we conclude that
$(Q^T)^*=Q^T$.
\smallskip

\noindent$\bullet$\,\,\,For $f \in {\cal M}^T \subset {\cal F}^T$, define
a function
 \begin{equation}\label{u^f new}
u^f(x,t)\,:=\,\left({\cal W}^T {\cal T}^T_{T-t} f\right)(x) \qquad
{\rm in} \,\,\,\overline{\Omega^T} \times [0,T]\,.
 \end{equation}
The definitions of the operators imply
 \begin{align*}
& \left[\Delta - Q^T\right]u^f(\cdot,t)=L^T u^f(\cdot,t)=L^T{\cal
W}^T {\cal T}^T_{T-t} f={\cal W}^T \left[{\cal
T}^T_{T-t} f\right]_{tt}=\\
& = [{\cal W}^T {\cal T}^T_{T-t} f]_{tt}\,=\,u^f_{tt}(\cdot,t)\,.
 \end{align*}
Thus, $u^f$ satisfies the equation
 \begin{equation}\label{E1}
u_{tt} - \Delta u + Q^T u\,=\,0 \qquad {\rm in} \,\,\,\Omega^T
\times (0,T)\,,
 \end{equation}
By (\ref{cal W^T new}) for $\sigma=\Gamma$, we have $\rm supp\,
u^f(\cdot,t)\subset \overline{\Omega^t}$, i.e., $u^f$ satisfies
the Cauchy condition
 \begin{equation}\label{E2}
u\big|_{t<\tau(x)}=0 \qquad {\rm in} \,\,\,\overline{\Omega^T}
\times [0,T]\,.
 \end{equation}
In the mean time, (\ref{boundary condition}) easily implies that
$u^f$ obeys
 \begin{equation}\label{E3}
\partial_\nu u=f \qquad {\rm on} \,\,\,\overline{\Sigma^T}\,.
 \end{equation}

\noindent$\bullet$\,\,\, Show that $Q^T$ is a multiplication by
bounded function. The proof follows the idea of \cite{BSobolev}.
 \begin{Lemma}\label{L3}
There is a (real) function $q \in L_\infty(\Omega^T)$ such that
$Q^T y= qy$ holds for $y \in {\cal H}^T$.
 \end{Lemma}
$\square$\,\,\,{1.}\,Choose a $\sigma \subset \Gamma$ and $f \in
{\cal F}^{T,\,\xi}_\sigma \cap {\cal M}^T$. By condition {\it 4} and (\ref{cal W^T
new}), we have $u^f(\cdot,T)\in {\cal H}^\xi_\sigma \cap H^2(\Omega^T)$. Hence,
$\Delta u^f(\cdot,T)\in {\cal H}^\xi_\sigma$. In the mean time, we have $f_{tt}
\in {\cal F}^{T,\,\xi}_\sigma \cap {\cal M}^T$ that implies
$u^f_{tt}=L^Tu^f(\cdot,T)={\cal W}^T f_{tt}\,\overset{(\ref{cal
W^T new})}\in\,{\cal H}^\xi_\sigma$. Therefore,
$Q^Tu^f(\cdot,T)\overset{(\ref{E1})}=\Delta
u^f(\cdot,T)-u^f_{tt}\in\,{\cal H}^\xi_\sigma$. Thus, $Q^T{\cal W}^T {\cal F}^{T,\,\xi}_\sigma
\subset {\cal H}^\xi_\sigma$ holds. Since ${\cal W}^T {\cal F}^{T,\,\xi}_\sigma$ is dense in
${\cal H}^\xi_\sigma$ (see (\ref{cal W^T new})), we conclude that $Q^T{\cal H}^\xi_\sigma
\subset {\cal H}^\xi_\sigma$. The latter leads to $Q^T[{\cal H}^T \ominus{\cal H}^\xi_\sigma] \subset
[{\cal H}^T \ominus{\cal H}^\xi_\sigma]$ by virtue of the symmetry $(Q^T)^*=Q^T$.
Hence, the subspaces ${\cal H}^\xi_\sigma$ {\it reduce} $Q^T$ that is equivalent
to the commutation
 \begin{equation}\label{commutation Q^TG=GQ^T}
Q^T G^\xi_\sigma\,=\,G^\xi_\sigma Q^T\,, \,\,\, \qquad \sigma
\subset \Gamma,\,\,\,\,0 \leqslant \xi \leqslant T\,,
 \end{equation}
where $G^\xi_\sigma$ projects in ${\cal H}^T$ onto ${\cal H}^\xi_\sigma$, i.e., cuts
off functions on $\Omega^\xi_\sigma$.
\smallskip

{2.}\,As is easy to verify, an operator $\tau^T_\sigma: {\cal H}^T
\to {\cal H}^T$,
 \begin{equation}\label{tau sigma def}
\tau^T_\sigma y := \left[\int_{[0,T]}\xi\,d
G^\xi_\sigma\right]y=\left[\lim \limits_{r_\Xi \to 0}\sum
\limits_{i=1}^N \xi_i\,[G^{\xi_i}_\sigma-G^{\xi_{i-1}}_\sigma]
\right]y
 \end{equation}
(the sums converge by the operator norm) acts by the rule
  $$
\tau^T_\sigma y= \begin{cases} {\rm d}(\cdot, \sigma)y  &
{\rm in}\,\,\,\Omega^T_\sigma\\
0 &
{\rm in}\,\,\,\Omega^T \setminus \Omega^T_\sigma\\
                      \end{cases}\,,
 $$
i.e., multiplies functions by the distance to $\sigma$ and, then,
cuts off on $\Omega^T_\sigma$ \cite{BSobolev}. As a consequence,
an operator
 $$
\hat \tau^T_\sigma\,:=\,\tau^T_\sigma y\,+\,T({\mathbb
I}_{{\cal H}^T}-G^T_\sigma)y
 $$
multiplies functions by the {\it continuous} function ${\rm
d}^T_\sigma(\cdot):=\max\{{\rm d}(\cdot, \sigma),T\}$. In the mean
time, (\ref{commutation Q^TG=GQ^T}) implies
 \begin{equation}\label{Q^T tau = tau Q^T}
Q^T\hat \tau^T_\sigma\,=\,\hat \tau^T_\sigma Q^T\,,\qquad \sigma
\subset \Gamma,\,\,\,\,0 \leqslant \xi \leqslant T\,,
 \end{equation}
because the sums in (\ref{tau sigma def}) do commute with all
$G^\xi_\sigma$.
\smallskip

{3.}\,\,Fix a (small) $\delta>0$. A simple geometric fact is that
the functions $\{{\rm d}^T_\sigma\,|\,\,\sigma \subset \Gamma\}$
separate points in $\Omega^{T-\delta}$ and vanish simultaneously
in no point $x_0 \in \overline{\Omega^{T-\delta}}$. Hence, a
family $\{{\rm d}^T_\sigma\,|\,\,\sigma \subset \Gamma,\,\,0
\leqslant \xi \leqslant T\}$ generates the continuous function
algebra $C(\overline{\Omega^{T-\delta}})$ \cite{BSobolev}.

Correspondingly, an operator family
$\{\hat\tau^T_\sigma\,|\,\,\sigma \subset \Gamma,\,\,0 \leqslant
\xi \leqslant T\}$ generates the operator (sub)algebra ${\mathfrak
C}(\overline{\Omega^{T-\delta}}) \subset {\mathfrak B}({\cal
H}^T)$ of multiplications by continuous functions. As a
consequence of (\ref{Q^T tau = tau Q^T}), we have $Q^T{\mathfrak
C}(\overline{\Omega^{T-\delta}})={\mathfrak
C}(\overline{\Omega^{T-\delta}})Q^T$ that is possible if and only
if $Q^T$ is also a multiplication by a function $q$.

Since $Q^T$ is bounded, we have $q\in
L_\infty(\Omega^{T-\delta})$. By arbitrariness of $\delta$, we get
$q\in L_\infty(\Omega^{T})$. \,\,\,$\blacksquare$
\medskip

\noindent$\bullet$\,\,\,With the above determined function $q$ one
associates the system $\alpha^T$ of the form
(\ref{B1})--(\ref{B3}). Such a system possesses its own operators
$W^T$ and $C^T$. Show that $W^T={\cal W}^T$ and $C^T={\cal C}^T$.

Since the problems (\ref{B1})--(\ref{B3}) and
(\ref{E1})--(\ref{E3}) (with $Q^T=q$) are identical and uniquely
solvable, their solutions (for the same $f$'s) coincide. Writing
the first relation of (\ref{Delay rel}) in the form
$u^f(\cdot,t)=W^T {\cal T}^T_{T-t} f$ and comparing with (\ref{u^f
new}), we see that $W^T={\cal W}^T$ holds.

By the latter equality and (\ref{{cal W}^*{cal W}= cal C}), we
have
 \begin{equation}\label{C^T=calC^T}
{\cal C}^T\,=\,({\cal W}^T)^*{\cal W}^T\,=\,(W^T)^*W^T\,=\,C^T\,.
 \end{equation}

\noindent$\bullet$\,\,\,System (\ref{E1})--(\ref{E3}) (with
$Q^T=q$) possesses the extended response operator $R^{2T}$. Here
we prove the equality $R^{2T}={\cal R}^{2T}$ that completes the
proof of the Theorem.

Begin with two lemmas of general character. The lemmas deal with a
Hilbert space ${\cal F}=L_2([0,2T];\,{\cal E})$ (with the Lebesgue
measure $dt$), where $\cal E$ is an auxiliary Hilbert space. By
${\cal F}_\pm$ we denote the subspaces of functions, which are
even and odd with respect to $t=T$. So, the decompositions ${\cal
F}\,=\,{\cal F}_+\oplus{\cal F}_-$ holds. Let
 $$
{\cal F}^{[a,b]}\,:=\,\{f \in {\cal F}\,|\,\,\rm supp\, f \subset
[a,b]\}\,, \qquad 0 \leqslant  a<b \leqslant 2T\,.
 $$
 \begin{Lemma}\label{L4}
If a bounded operator $N: {\cal F}\to \cal F$ satisfies
 \begin{equation}\label{N conditions}
N {\cal F}_\pm \,\subset\, {\cal F}_\pm\,; \quad N{\cal
F}^{[a,2T]}\,\subset\,{\cal F}^{[a,2T]}\,, \,\,\,0  \leqslant a \leqslant 2T
 \end{equation}
then it is local, i.e., preserves the support of functions:
 \begin{equation}\label{N local}
N {\cal F}^{[a,b]}\,\subset\,{\cal F}^{[a,b]}\,,\qquad 0 \leqslant  a<b
\leqslant 2T\,.
 \end{equation}
 \end{Lemma}
$\square$\,\,\,\,\,\,\,{\bf 1.}\,\,Representing ${\cal F}={\cal
F}^{[0,T]}\oplus{\cal F}^{[T,2T]}$ and $f=f_1+f_2$ with $f_1
\in{\cal F}^{[0,T]},\,f_2 \in{\cal F}^{[T,2T]}$, we identify
$f\equiv\langle f_1,f_2\rangle$.

Introduce an isometry $Y: {\cal F}^{[0,T]}\to{\cal F}^{[T,2T]}$ by
 $$
(Yf)(t)\,:=\,f(2T-t)\,, \qquad T \leqslant t \leqslant 2T\,.
 $$
Obviously, one has ${\cal F}_\pm =\{\langle f,\pm
Yf\rangle\}\,|\,\,f \in {\cal F}^{[0,T]}\}$. Since $N$ preserves
the evenness/oddness, there are two operators $k,l: {\cal
F}^{[0,T]}\to{\cal F}^{[0,T]}$ such that
 \begin{equation}\label{Aux 1}
N\langle f,Yf\rangle = \langle kf, Ykf\rangle \qquad {\rm and}
\qquad N\langle f,-Yf\rangle = \langle lf, -Ylf\rangle\,.
 \end{equation}
Show that $k=l$. For a $g \in {\cal F}^{[0,T]}$, one has
 \begin{align}
\notag & 2N\langle 0, Yg\rangle = N\left[ {\langle g, Yg\rangle -
\langle g, -Yg\rangle}\right] \overset{(\ref{Aux 1})}= \langle kg,
Ykg\rangle - \langle lg,
-Ylg\rangle=\\
\label{Aux 2}& = \langle [k-l]g, Y[k+l]g\rangle\,.
 \end{align}
In the mean time, we have $\langle 0,Yg\rangle \in {\cal
F}^{[T,2T]}$ and, hence, $N\langle 0,Yg\rangle \in {\cal
F}^{[T,2T]}$ holds by (\ref{N conditions}). By the latter,
$2N\langle 0, Yg\rangle$ must be of the form $\langle 0,
...\rangle$, i.e., $[k-l]g=0$ is valid and implies $k=l=:m$.
\smallskip

\noindent{\bf 2.}\,\,\,Putting $g=Y^{-1}h$ in (\ref{Aux 2}), we
get
 \begin{equation}\label{Aux 3}
N\langle 0, h\rangle\,=\,\langle 0,YmY^{-1}h\rangle\,.
 \end{equation}
In the mean time, we have
 \begin{align*}
\notag & 2N\langle g, 0\rangle = N\left[ {\langle g, Yg\rangle +
\langle g, -Yg\rangle}\right] \overset{(\ref{Aux 1})}= \langle mg,
Yg\rangle + \langle mg, -Ymg\rangle = \\
 & =\,2\langle mg, 0\rangle\,.
 \end{align*}
Combining the latter with (\ref{Aux 3}), we arrive at the
representation
 \begin{equation}\label{Aux 4}
N\langle g, h\rangle\,=\,\langle mg,YmY^{-1}h\rangle\,.
 \end{equation}

\noindent{\bf 3.}\,\,\,Such a representation easily provides the
following fact: operator $N$ acts locally in $[0,2T]$ if and only
if operator $m$ is local in  $[0,T]$. Show that the latter does
occur.

Let ${\rm supp\,}f \subset [a,b]\subset [0,T]$, so that
 $
f\big|_{0\leqslant t<a}=0 \qquad {\rm and} \qquad
f\big|_{b<t\leqslant 2T}=0
 $
holds. The first equality means that $f \in {\cal F}^{[a,2T]}$,
implies $Nf \in {\cal F}^{[a,2T]}$ by (\ref{N conditions}) and,
thus, provides $Nf\big|_{0\leqslant t<a}=0$. Hence, with regard to
$f \equiv \langle f,0\rangle$, we have
 \begin{align*}
0=Nf\big|_{0\leqslant t<a}\equiv \left[N\langle
f,0\rangle\right]\big|_{0\leqslant t<a}\overset{(\ref{Aux
4})}=\langle mf,0\rangle\big|_{0\leqslant t<a}\equiv
mf\big|_{0\leqslant t<a}\,,
 \end{align*}
i.e., $m$ does not extend support to the left.

By the choice of $f$, one has ${\rm supp\,}Yf \subset[2T-b,
2T-a]$, so that $Yf \in {\cal F}^{[2T-b,2T]}$. The latter implies
$NYf \in {\cal F}^{[2T-b,2T]}$ in accordance with (\ref{N
conditions}). Hence, we have
 \begin{align*}
& 0=NYf\big|_{0\leqslant t<2T-b}\equiv \left[N\langle
0,Yf\rangle\right]\big|_{0\leqslant t<2T-b}\overset{(\ref{Aux
4})}=\langle 0,Ymf\rangle\big|_{0\leqslant t<2T-b}\equiv\\
& \equiv\,Ymf\big|_{0\leqslant t<2T-b}\,.
 \end{align*}
Therefore, $mf\big|_{t>2T-b}=0$, i.e., $m$ does not extend support
to the right. Thus, $m$ acts locally and, eventually, $N$ is
local.\,\,\,\,\,$\blacksquare$

In fact, the boundedness of $N$ is not substantial and the proof
(mutatis mutandis) is available for a wider class of operators.
 \begin{Lemma}\label{L5}
If an operator $N$ satisfies (\ref{N conditions}) and is compact
then $N=\mathbb O$.
 \end{Lemma}
$\square$\,\,\,\,A projection $X^{[a,b]}$ in $\cal F$ onto ${\cal
F}^{[a,b]}$ cuts off functions on $[a,b]$. The complement
projection $X^{[a,b]}_\bot={\mathbb I}-X^{[a,b]}$ cuts off on
$[0,a]\cup [b,2T]$. By Lemma \ref{L4}, we have
 $$
NX^{[a,b]}=X^{[a,b]}NX^{[a,b]} \quad {\rm and} \quad
NX^{[a,b]}_\bot=X^{[a,b]}_\bot NX^{[a,b]}_\bot\,.
 $$
Summing up, we get
 $
N=X^{[a,b]}NX^{[a,b]}+X^{[a,b]}_\bot NX^{[a,b]}_\bot
 $
that leads to
 $$
NX^{[a,b]}=X^{[a,b]}N\,, \qquad N^*X^{[a,b]}=X^{[a,b]}N^*
 $$
and, eventually, implies
 \begin{equation}\label{commutation N*NX=XN*N}
N^*NX^{[a,b]}\,=\,X^{[a,b]}N^*N\,.
 \end{equation}
In the mean time, operator $N^*N$ is self-adjoint and compact. Let
$\lambda \in \mathbb R$ be its eigenvalue, ${\cal D}_\lambda$ the
corresponding eigensubspace. By (\ref{commutation N*NX=XN*N}), we
have $X^{[a,b]}{\cal D}_\lambda\subset{\cal D}_\lambda$ that leads
to ${\rm dim\,}{\cal D}_\lambda=\infty$. The latter is possible
only for ${\cal D}_0={\rm Ker\,}N^*N$. Thus, the spectrum of
$N^*N$ is exhausted by $\lambda=0$. Hence, $N^*N=\mathbb O$.
Therefore, $N=\mathbb O$. \,\,\,\,\,$\blacksquare$

\noindent$\bullet$\,\,\,Now, we are ready to complete the proof of
Theorem \ref{T1}. Return to our system (\ref{E1})--(\ref{E3})
(with $Q^T=q$). Recall that $S^T: {\cal F}^T \to {\cal F}^{2T}$
extends controls from $[0,T]$ to $[0,2T]$ by oddness with respect
to $t=T$. We regard ${\cal F}^{2T}=L_2(\Sigma^{2T})$ as the space
$L_2([0,2T];{\cal E})$ with ${\cal E}=L_2(\Gamma)$. Let ${\cal
F}^{2T}_\pm$ be the subspaces of the even and odd functions, so
that the decomposition
 \begin{equation*}
{\cal F}^{2T}\,=\,{\cal F}^{2T}_+\oplus{\cal F}^{2T}_-
 \end{equation*}
occurs. The embedding $J^{2T}{\cal F}^{2T}_-\subset {\cal
F}^{2T}_+$ holds and is dense. Also, one has $Y^{2T}{\cal
F}^{2T}_\pm={\cal F}^{2T}_\pm$.

Denote $N:=R^{2T}-{\cal R}^{2T}$ With regard to (\ref{C^T via
R^2T}) and (\ref{cal C^T}), the equality (\ref{C^T=calC^T}) leads
to
 $$
(NJ^{2T}S^T f, S^T g)_{{\cal F}^{2T}}=0
 $$
for all $f,g \in {\cal F}^T$. It shows that the embedding
 \begin{equation*}
N{\cal F}^{2T}_+\,\subset\,{\cal F}^{2T}_+
 \end{equation*}
holds and evidently implies $Y^{2T}N{\cal
F}^{2T}_+\,\subset\,{\cal F}^{2T}_+$. In the mean time, operator
$Y^{2T}N$ is self-adjoint: see (\ref{R^2T rel}) and (\ref{cal R}).
Therefore, it is reduced by the even/odd subspaces: $Y^{2T}N{\cal
F}^{2T}_\pm\,\subset\,{\cal F}^{2T}_\pm$. The latter leads to
 \begin{equation}\label{Embedding 1}
N{\cal F}^{2T}_\pm\,\subset\,{\cal F}^{2T}_\pm\,.
 \end{equation}

On the other hand, the shift invariance (\ref{R^2T rel}) and
(\ref{cal R}) implies
 \begin{equation}\label{Embedding 2}
N{\cal F}^{2T,\,\xi}\,\subset{\cal F}^{2T,\,\xi}\,, \qquad 0
\leqslant \xi \leqslant 2T\,.
 \end{equation}
Joining (\ref{Embedding 1}) with (\ref{Embedding 2}) and applying
Lemma \ref{L5}, we arrive at $N=\mathbb O$ that is $R^{2T}={\cal
 R}^{2T}$. Theorem \ref{T1} is proved.\,\,\,\,$\blacksquare$

\section{Comments, doubts, philosophy}\label{Comments}
$\bullet$\,\,\,A characterization of data for an inverse problem
is a list of conditions providing its solvability. The reasonable
requirement to any characterization is to be checkable and
possibly simple. As we guess, the only reasonable understanding of
`a condition is checkable' is that it can be verified {\it before}
(without) solving the inverse problem. Formally, the conditions
{\it 1--7} of Theorem \ref{T1} satisfy such a requirement because
they do not use the knowledge of the potential $q$. However,
comparing these conditions with the procedure {\it Step 1--4}, it
is easy to recognize that to check {\it 1--7} is almost the same
as to recover $q$. Conditions {\it 1--7} just provide the
procedure to be realizable. In such a situation, can one claim
that {\it 1--7} is an efficient characterization?

And what is `efficient'? For instance, the key step of the
procedure, as well as the characterization, is constructing the
operator integral (\ref{cal U^T}). If it is at our disposal, we
get $W^T$, recover the waves $u^f$, and are able to check {\it
5--7}. In the mean time, having $u^f$ one doesn't need to check
anything more but can just determine $q$ from the wave equation.
So, can one regard the required in {\it 3} convergence as an
efficiently checkable condition? We don't have a convincible
answer.

Also, can one avoid so long list of conditions and invent
something simpler and better? \footnote{Actually, a long list of
the characterization conditions is not something unusual: see,
e.g., the conditions on a spectral triple corresponding to a
Riemannian manifold in \cite{Connes}.} We are rather sceptical and
the following is some reasons for scepticism.
\smallskip

$\bullet$\,\,\,The evolution of system (\ref{B1})--(\ref{B3}) is
governed by the operator $L_q=-\Delta+q$ and Neumann controls
$f=\partial_\nu u\big|_{\Sigma^T}$. Both of them are of very
specific type. We mean, replacing them by
$L_Q=-\sum_{i,j}\partial_{x_i}a^{ij}\partial_{x^j}+Q$ (with
possibly nonlocal and time dependent $Q$) and, let say,
$f=[\partial_\nu u+\varkappa u]\big|_{\Sigma^T}$, we'd got a
system with the data $R^{2T}_Q$ of the properties quite analogous
to $R^{2T}_q$. Therefore, the data characterization has to select
$R^{2T}_q$ from a large reserve of the response operators
$R^{2T}_Q$. It is such a selection, which the conditions {\it
1--7} do implement. Namely, the selection works as follows.
\smallskip

\noindent$\star$\,\,\,Conditions {\it 1,\,2} appear at very
general level of an abstract dynamical system with boundary
control (DSBC) associated with a time-independent boundary triple
\cite{DSBC}. Such a system necessarily satisfies  (\ref{cal R})
and (\ref{cal C^T}).
\smallskip

\noindent$\star$\,\,\, In {\it 3}, convergence of the operator
integral to an isometric operator is a specific feature of {\it
hyperbolic} DSBC's obeying the finiteness of domain of influence
principle. System $\alpha^T$, which we deal with, is hyperbolic,
and the characterization must provide such a property.

Also, as was noticed in sections \ref{subsec System alpha^T},
\ref{subsec Devices} (see (\ref{C^T triang factor}), (\ref{V^T
def})), the amplitude integral is connected with a triangular
factorization. One of the form of the classical factorization
problem is to recover a triangular operator via its imaginary
(anti-Hermitian) part. It is solved by the use of the so-called
triangular truncation transformer \cite{GohKr}, which is a kind of
an operator integral. Its convergence provides a solvability {\it
criterium} to the factorization problem for a class of Fredholm
operators \cite{GohKr}.

So, imposing condition {\it 3}, we follow the classicists. By the
way, our construction (\ref{AI-2}) is available for a wider class
of operators \cite{BPush}.
\smallskip

\noindent$\star$\,\,\,The characterization should specify a
regularity class of potentials, which we deal with. Condition {\it
4}, roughly speaking, rejects strongly singular potentials.
\smallskip

\noindent$\star$\,\,\,Condition {\it 5} excludes another types of
boundary conditions like $f=[\partial_\nu u+\varkappa
u]\big|_{\Sigma^T}$. The Neumann condition is rather specific. In
contrast to the Dirichlet condition, which is connected with a
Friedrichs operator extension, the Neumann one is not of invariant
meaning. The characterization has to take this fact into account.
Perhaps, one can specify the boundary condition right from
$R^{2T}$, without constructing $W^T$. It would be welcome.
\smallskip

\noindent$\star$\,\,\,A discussable question is whether condition
{\it 6} may be efficiently checked. However, (\ref{cal W^T new})
is also unavoidable: it is the condition, which provides a {\it
locality} of the potential.
\smallskip

\noindent$\star$\,\,\,Assume for a while that $q \in
L_2(\Omega)\setminus L_\infty(\Omega)$, so that the multiplication
by $q$ is an {\it unbounded} operator. However, system $\alpha^T$
with such a potential does possess all the properties specified by
conditions {\it 1--6}. In the mean time, the characterization must
reject such a case. We see no option to do it except of imposing
(\ref{Estimate Q^T}).
\smallskip

So, all the conditions {\it 1--7} are independent and, therefore,
unavoidable. We are forced to accept so long list of conditions
just because we deal with a very specific class of dynamical
systems. The more specific is the class, the more words is
required for its description. The converse is also true: to be the
response operator of an abstract DSBC, it suffices for ${\cal
R}^{2T}$ to satisfy nothing but (\ref{cal R}) and (\ref{cal C^T})
\cite{DSBC}.
\medskip

$\bullet$\,\,\,A determination of $q$ from $R^{2T}$ is
conventionally regarded as an {\it over-determined} problem. The
reason is the following. One can represent
 $$
\left(R^{2T} f\right)(\gamma,t)=\int_{\Sigma^t}r(t-s, \gamma,
\gamma')\,f(\gamma',s)\,d\Gamma_{\gamma'}\,ds
 $$
with a (generalized) kernel $r(t, \gamma, \gamma')$. The
convolution form with respect to time is  a consequence of the
shift invariance (\ref{R^2T rel}). Bearing in mind that
$\gamma=\{\gamma^1, \gamma^2, \dots, \gamma^{n-1}\}$, one regards
$r$ as a function of $1+2(n-1)=2n-1$ variables, whereas a {\it
local} potential $q=q(x^1, x^2, \dots, x^{n})$ depends on $n$
variables only. Thus, for $n \geqslant 2$ the data array is of
higher dimension than the array of parameters under determination
`that is not natural'.

Actually, on our opinion, in {\it multidimensional} problems such
a counting parameters is not quite relevant and reliable
\footnote{for instance, how to count the parameters if we need to
recover from $R^{2T}$ not a function (potential) but a Riemannian
manifold, as in \cite{BIP07}?}. Nevertheless, the question arises:
Does the characterization {\it 1--7} `kill' unnecessary parameters
and, if yes, in which way? The possible answer is the following.

There is a sharp  {\it necessary} condition related with a
locality of potential. Let $\tilde {\cal P}^{T,\,\xi}_\sigma$ be
the projection in ${\cal F}^T$ onto the subspace $\overline{[{\cal
C}^T]^{\frac{1}{2}}{\cal F}^{T,\,\xi}_\sigma}$. Such a projection
is unitarily equivalent (via the isometry $(I^T)^*{\cal A}^T$: see
(\ref{Basic for cal W^T})) to the projection onto $\overline{{\cal
W}^T{\cal F}^{T,\,\xi}_\sigma}$. By (\ref{cal W^T new}), the
latter projection coincides with the `geometric' projection
$G^\xi_\sigma$, which cuts off functions onto $\Omega^\xi_\sigma$.
The geometric projections for all $\sigma$ and $\xi$ commute. As a
result, we arrive at the following condition: {\it the projection
family $\{\tilde {\cal P}^{T,\,\xi}_\sigma\,|\,\,\sigma \subset
\Gamma,\,0 \leqslant \xi \leqslant T\}$ must be commutative}.
Analyzing the proof of Theorem \ref{T1}, we see that it is the
condition, which forces the `potential' $Q$ to be a multiplication
by $q$ and, thus, rejects unnecessary variables. However, the
rejection mechanism is not well understood yet and we hope to
clarify it in future.

\end{document}